\documentclass[a4paper,12pt]{article}
\usepackage[english]{babel}
\usepackage[latin1]{inputenc}
\usepackage[T1]{fontenc}
\usepackage{amsmath}
\usepackage{amsfonts}
\usepackage{dsfont}
\usepackage{amsthm}
\usepackage{amssymb}
\title{On the measures of large entropy on a positive closed current}
\usepackage{dsfont}
\author{Henry de Th\'elin and Gabriel Vigny}
\begin{document}
\newtheorem{Theorem}{Theorem}
\newtheorem*{Theorem*}{Theorem}
\newtheorem{theorem}{Theorem}[section]
\newtheorem{proposition}[theorem]{Proposition}
\newtheorem{defi}[theorem]{Definition}
\newtheorem{corollary}[theorem]{Corollary}
\newtheorem{Hypothesis}[theorem]{Hypothesis}
\newtheorem{lemma}[theorem]{Lemma}
\newtheorem{Remark}[theorem]{Remark}
\newcommand{\U}{\mathcal{U}}
\newcommand{\C}{\mathcal{C}}
\renewcommand{\P}{\mathbb{P}}
\renewcommand{\phi}{\varphi}
\newcommand{\Cc}{\mathbb{C}}
\newcommand{\Nn}{\mathbb{N}}
\newcommand{\Rr}{\mathbb{R}}
\newcommand{\Qq}{\mathbb{Q}}
\newcommand{\Zz}{\mathbb{Z}}
\newcommand{\Acal}{\mathcal{A}}
\newcommand{\Bcal}{\mathcal{B}}
\newcommand{\Dcal}{\mathcal{D}}
\newcommand{\Ecal}{\mathcal{E}}
\newcommand{\Gcal}{\mathcal{G}}
\newcommand{\Hcal}{\mathcal{H}}
\newcommand{\Lcal}{\mathcal{L}}
\newcommand{\Ical}{\mathcal{I}}
\newcommand{\Pcal}{\mathcal{P}}
\newcommand{\Qcal}{\mathcal{Q}}
\newcommand{\Scal}{\mathcal{S}}
\newcommand{\Zcal}{\mathcal{Z}}
\hyphenation{plu-ri-sub-har-mo-nic}
\date{}

\maketitle
\begin{abstract}  

Let $f:X\to X$ be a dominating meromorphic map of a compact Kähler surface of large topological degree.  Let $S$ be a positive closed current on $X$ of bidegree $(1,1)$. We consider an ergodic measure $\nu$ of large entropy supported by $\mathrm{supp}(S)$. Defining dimensions for $\nu$ and $S$, we give inequalities à la Mañé involving the Lyapunov exponents of $\nu$ and those dimensions. We give dynamical applications of those inequalities.

 \end{abstract}


\section{Introduction}
Let $X$ be a compact Kähler surface endowed with a Kähler form $\omega$. Let $f$ be a dominating meromorphic map of $X$. We define the quantity:
$$ \lambda_1(f^n):=\int_X (f^n)^* (\omega) \wedge \omega = \int_X (f^n)_* (\omega) \wedge \omega .$$
 Recall that the first dynamical degree $d_1$ is the well defined quantity: 
$$d_1:= \lim_n \left(\lambda_1(f^n) \right)^{1/n}.$$
Assume that the topological degree $d_t$ satisfies $d_t >d_1$ (e.g. if $f$ is a holomorphic map of $\P^2$ of algebraic degree $d$, then $d_t= d^2> d=d_1$). We let $I$ denote the indeterminacy set of $f$, $\C$ its critical set and $\Acal:= I \cup \C$. Let $S$ be a positive closed current on $X$ of bidegree $(1,1)$, we normalize it so that its mass satisfies $\|S\|=1$. The purpose of the article is to study the measures of large entropy whose support is contained in  $\mathrm{supp}(S)$. Let $\Lambda$ be a subset of $\mathrm{supp}(S)$.
In what follows, we define
\begin{equation*}
 \bar{d}_S(x) :=\overline{\lim_{r\to 0}} \frac{\log (S\wedge \omega (B(x,r)))}{\log r}  
\end{equation*}
and we let $\bar{d}_S$ be the supremum of the $\bar{d}_S(x)$, $x\in \Lambda$, in particular:
\begin{equation}\label{hypotheseS}
\forall x \in \Lambda, \ \overline{\lim_{r\to 0}} \frac{\log (S\wedge \omega (B(x,r)))}{\log r}\leq \bar{d}_S. 
\end{equation}
Let $\nu$ be an $f$-invariant ergodic measure. By \cite{Dup, DT1}, we have that if the entropy of $\nu$ satisfies $h_{\nu}(f) > \log d_1$ (that is what we mean by large entropy) and $\log d(x, I) \in L^1(\nu)$ then the Lyapunov exponents are well defined and satisfy  $\chi_1 \geq \chi_2>0$. The following theorems are the main results of the article:
\begin{Theorem}\label{Theorem1} Let $\nu$ be an $f$-invariant ergodic measure with $\nu(\Lambda)=1$ such that $\log d(x, I) \in L^1(\nu)$ and $h_{\nu}(f) > \log d_1$. Then we have the inequality:
\begin{equation}\label{equation1}
\chi_1\bar{d}_S -2 \chi_2 \geq h_{\nu}(f) -\log d_1.
\end{equation}
\end{Theorem}
\begin{Theorem}\label{Theorem2} We keep the hypotheses and notations of Theorem \ref{Theorem1}. Let $\underline{d}_{\nu}$ be such that for $\nu$-almost every $x$ in $\Lambda$:
$$\liminf_{r\to 0} \frac{\log(\nu(B(x,r))) }{\log r } \geq \underline{d}_{\nu}.$$  
Then:
\begin{equation}\label{equation2}
\chi_1(\bar{d}_S-\underline{d}_{\nu}) \geq 2 \chi_2-\log d_1.
\end{equation}
\end{Theorem}
By \cite{Dup}[Corollary 2], we know that, for an endomorphism $f$ of $\P^2$, there exists some $\varepsilon>0$ such that for any $h\in [\log d_t-\varepsilon, \log d_t]$ there exists a measure of entropy $h$ with support contained in the support of the Green current $T$. So the above theorems can be applied to many examples. 
Theorem \ref{Theorem2} is actually finer than Theorem \ref{Theorem1} and its proof is given by a modification of the one of Theorem  \ref{Theorem1}. Nevertheless both statements have their uses. We start by giving applications of the previous theorems then we sketch the proofs. We recall some fact on Pesin's theory and give a theorem of graph transform. Finally, we give the proof of Theorem \ref{Theorem1} and then of Theorem \ref{Theorem2}.

\section{Applications}

 \noindent {\bf $\bullet$ Sharpness and lower bound of the higher Lyapunov exponent.} 

Observe first that the inequalities (\ref{equation1}) and (\ref{equation2}) are sharp. Let $g$ be a rational map of $\P^1$ of degree $d$ and let $\nu_1$ be the measure of maximal entropy $\log d$. Mané's results (\cite{Man}) imply that $\bar{d}_{\nu_1}\chi=\underline{d}_{\nu_1} \chi= \log d $ where the Lyapunov exponent $\chi >0$ is well defined.
 Consider $f=(g,g)$ acting on $(\P^1)^2$. Take $S := (\pi_1)^*(\nu_1)$ ($\pi_1$ is the projection on the first coordinate) and $\nu := \nu_1 \otimes \nu_1$. Take $\Lambda=\pi_1^{-1}(R)$ where $R$ is the set of $\P^1$ of full measure for $\nu_1$ such that $\bar{d}_{\nu_1}(x)=\underline{d}_{\nu_1}(x)$ for $x\in R$.
Then $\bar{d}_S= \underline{d}_{\nu_1} +2$, $\chi_1=\chi_2 = \chi$,  
$h_{\nu}(f)=2  \log d$ and $d_1=d$ and one sees that (\ref{equation1}) and (\ref{equation2}) are equalities here. \\
 
Finally, we know that $2\chi_2 \geq h_\nu - \log d_1$ (\cite{DT1}[Corollaire 2]). In particular, Theorem \ref{Theorem1} implies here (which is interesting as soon as $\bar{d}_S<4$):
\begin{corollary}
With the above hypothesis, we have $\chi_1 \geq \frac{2}{\bar{d}_S}(h_\nu - \log d_1)$. 
\end{corollary}

\noindent {\bf $\bullet$ A necessary condition for the lower Lyapunov exponent to be minimal} 

We give some corollaries of Theorem \ref{Theorem2}. 
\begin{corollary}  If $\bar{d}_S = \underline{d}_{\nu}$ then $ \chi_2 \leq \frac{\log d_1 }{2}.$
\end{corollary}
\begin{corollary}\label{above}  If $f$ is holomorphic in $X = \P^2$, $\nu= \mu$ is the measure of maximal entropy and if $\bar{d}_S=\underline{d}_{\nu}$ then $\chi_2= \frac{\log d }{2}$.
\end{corollary}
\noindent {\it Proof of Corollary \ref{above}.} The inequality $\chi_1 \geq \chi_2 \geq \frac{\log d }{2}$ always stands by \cite{BrDu}, we conclude using the previous corollary.  \hfill $\Box$ \hfill \\

The hypotheses of Corollary \ref{above} are, in particular, satisfied  for Lattès maps (taking $S$ equal to the Green current) and actually, we know that for (and only for) Lattès maps, both Lyapunov exponents are minimal equal to  $ \frac{\log d }{2}$ (see \cite{BeDu}[Corollaire 1]).  \\

\noindent {\bf $\bullet$ Lower bound on $d_T$}

Consider again the case where $f$ is a holomorphic endomorphism of $\P^2$ of degree $d$. Take $S=T$ equal to the Green current and $\nu=\mu$ the measure of maximal entropy. Then Dupont proved in \cite{Dup2} that
$$ \underline{d}_\nu \geq   \frac{\log d }{\chi_1} + \frac{\log d }{\chi_2}  $$
(the equality was conjectured in \cite{BDM}). In particular, we get from Theorem \ref{Theorem2} that:
$$ \bar{d}_T \geq 2 \frac{\chi_2}{\chi_1} + \frac{\log d}{\chi_2} . $$
 Furthemore, minimizing the function $x\mapsto 2x/ \chi_1 + (\log d)/ x$ implies the bound:
$$ \bar{d}_T \geq 2 \sqrt{\frac{2 \log d }{\chi_1}}, $$
which only uses the largest Lyapunov exponent. \\

\noindent {\bf $\bullet$ Upper bound on the Hölder exponent of the Green function} 

Consider again the case where $f$ is holomorphic in $\P^2$ (for the sake of simplicity), $\nu= \mu$ is the measure of maximal entropy and $S=T$ is the Green current of $f$. Then it is known that $\mu = T \wedge T$ and $T$ has $\alpha$-Hölder continuous potentials (\cite{DS}). In particular, we take for $\Lambda$ the set of $x$ such that:
 $  \underline{d}_\mu(x) \geq \underline{d}_\mu  \ \mathrm{and} \ \bar{d}_\mu(x) \leq \bar{d}_\mu $, $x$ is generic for Pesin's theory and satisfies Brin-Katok formula. Then take $r$ small so that $\mu(B(x,r)) \geq r^{\overline{d}_\mu + \gamma}$ for some $\gamma$ arbitrarily small ($r$ just have to be small enough). In $B(x,2r)$, we write $T=dd^c \varphi$ where $\varphi$ is a $C^\alpha$-psh map. Let $0\leq \chi_r \leq 1$ be a smooth cut-off function equal to $1$ in $B(x,r)$ and $0$ outside  $B(x,2r)$, then $\|dd^c \chi_r \|_\infty \leq Cr^{-2}$ where $C$ does not depend on $r$. Then, by Stokes' formula:
\begin{align*}
r^{\bar{d}_\mu + \gamma} &\leq \int_{B(x,r)} T\wedge T \leq \int \chi_r T\wedge dd^c (\varphi-\varphi(x)) \leq  \int  (\varphi-\varphi(x)) T\wedge dd^c \chi_r \\
                               &\leq C'r^{\alpha -2} \int_{B(x,2 r)}   T\wedge \omega \leq C' r^{\alpha-2} (2r)^{ \bar{d}_T(x)-\gamma} 
\end{align*}
for $r$ small enough well chosen and for $C'$ a constant independent of $x$ (here we take an appropriate subsequence of $r$ going to $0$). Letting $r \to 0$ and $\gamma \to 0$ gives $\bar{d}_\mu \geq \bar{d}_T(x)+\alpha-2$ and taking the supremum for all $x$:
\begin{align}\label{dtdmu1}
  \bar{d}_\mu \geq \bar{d}_T+\alpha-2. 
\end{align}

In the same manner, we have by Brin-Katok's formula that $2\log d= h_\mu(f) =\lim_{r\to 0} \liminf_{n\to \infty} - \log(\mu(B_n(x,r)))/n$ where $B_n(x,r)$ is the Bowen ball and $x$ is $\mu$-generic. It is classical in that setting that (see \cite{HZ}) for $x$ generic, $\gamma$ arbitrarily small and $n$ large enough:
$$ B_n(x,r) \subset B(x,e^{-\chi_2 n+\gamma n})   . $$  
Arguing as above gives that ($C$ is a constant that does not depend on $x$ nor $n$):
\begin{align*}   
\mu(B_n(x,r)) & \leq  \mu(B(x,e^{-\chi_2 n+\gamma n}))  \\
&\leq  C (e^{-\chi_2 n+\gamma n})^{\alpha-2} (2e^{-\chi_2 n+\gamma n})^{ \bar{d}_T(x)-\gamma}, 
\end{align*} 
here, one needs to choose a suitable $n$ to get the $\limsup$. Taking the logarithm, dividing by $n$, letting $n\to \infty$, $\gamma \to 0$ and taking the supremum for all $x$ give:
\begin{align}\label{dtdmu2}
 2\log d \geq \chi_2 (\bar{d}_T +\alpha -2).
\end{align}
Recall that the measure $\mu$ is exact dimensional if $ \bar{d}_\mu= \underline{d}_\mu$ (this should hold in our case but it is unknown). Combining (\ref{dtdmu2}) with Theorem \ref{Theorem1} and (\ref{dtdmu1}) with  Theorem \ref{Theorem2}, we get the following corollary that gives bounds for $\alpha$ in term of the Lyapunov exponents and $\log d$ :
\begin{corollary}
With the above notations, one has:
$$\log d \left(\frac{2}{\chi_2} - \frac{1}{\chi_1}\right) + 2 \left(1 - \frac{\chi_2}{\chi_1}\right)  \geq \alpha . $$
Assume furthermore that $\mu$ is exact dimensional, then one has the sharper bound:
$$ \frac{\log d}{\chi_1} + 2 \left(1 - \frac{\chi_2}{\chi_1}\right)  \geq \alpha .$$ 
\end{corollary}
Observe that \cite{DS}[Proposition 1.18] gives a lower bound for $\alpha$ in the sense that for any $\alpha < \mathrm{min}(1, \log d/ \log (\lim \|Df^n\|_\infty^{1/n}))$ then $f$ is $\alpha$-Hölder, our result is a sort of reverse. \\

\noindent {\bf $\bullet$ Application to real measure for real rational map.} 

Consider a rational map $f_\Cc$ of $\P^2$ with real coefficient with $d_t > d_1$. In particular, it defines by restriction a rational map $f_\Rr$ of $\P^2(\Rr)$. Assume that $f_\Cc$ admits a measure $\nu$ with support in $\P^2(\Rr)$, such that $h_\nu(f_\Cc)> \log d_1$ and  $\log d(x, I) \in L^1(\nu)$. One can find a positive closed  current $S$ such that $\P^2(\Rr) \subset \mathrm{supp}(S)$ and $\bar{d}_S=3$ for $\Lambda = \P^2(\Rr)$ (take some fibration by complex lines above a $\P^1(\Rr)$). In particular, Theorem \ref{Theorem1} implies here:
\begin{corollary}
With the above hypothesis, one has $3\chi_1-2\chi_2 \geq h_\nu(f_\Cc)-\log d_1 $.
\end{corollary} 
Let $g$ be a Chebyshev's polynomial of degree $d$ and take $f=(g,g)$. Then $\mu_f=\mu_g\otimes \mu_g$ (where $\mu_g$, the measure of maximal entropy $\log(d)$ of $g$, has support in $\Rr$) has support in $\Rr^2$ and is of entropy $2\log d$. Recall that here $\chi_1=\chi_2=\log d$ hence we have equality in the above corollary. In particular, this shows that $\bar{d}_S=3$ is the best value one can achieve: there is no positive closed $(1,1)$ current $S$ in $\P^2(\Cc)$ such that $\P^2(\Rr) \subset \mathrm{supp}(S)$ and $\bar{d}_S<3$ for $\Lambda = \P^2(\Rr)$. \\

\noindent {\bf $\bullet$ Sets of small entropy.}  

Finally, we have the corollary which immediately follows from Theorem \ref{Theorem1}:
\begin{corollary}
Let $f$ be holomorphic on $X$ and $S$ be a current such that $\bar{d}_S=2$ for $\Lambda=\mathrm{supp}(S)$. Then there exists no invariant measure supported in $\mathrm{supp}(S)$ of entropy $>\log d_1$ and Lyapunov exponents satisfying $\chi_1=\chi_2$. 
\end{corollary}
It would be nice to remove the condition $\chi_1=\chi_2$ from the previous corollary. In fact, this was one of the initial motivations of this work: in \cite{Buf}, the author constructs a birational map of $\P^2$ such that the set where the dynamics may have entropy is contained in a pluripolar current (in a loose sense). Such examples could be counter examples to the entropy conjecture for birational map (existence of a measure of maximal entropy equal to $\log d$). A first step would be to understand  the holomorphic case. Hence we have the following questions:\\

\noindent {\bf Questions} 
\begin{itemize}
\item Let $f$ be a holomorphic map of $X$ and let $S$ be a $(1,1)$ positive closed current of Hausdorff dimension $2$, is there an invariant measure supported in $\mathrm{supp}(S)$ of entropy $>\log d_1$? 
\item  Let $f$ be a holomorphic map of $X$ and let $P$ be a compact pluripolar set of $X$, is there an invariant measure supported in $P$ of entropy $>\log d_1$? 
\end{itemize}
For the first question, the Hausdorff dimension could mean several things: the support of $S$ is of Hausdorff dimension
$2$ or  the measure $S\wedge \omega$ is of Hausdorff dimension $2$ in the sense that $\bar{d}_{S\wedge \omega}=2$ or 
 that $\underline{d}_{S\wedge \omega}=2$ or that $\mathrm{dim}_{\Hcal}(S\wedge \omega)=2$.  We expect that the answers are no (this is the case in dimension $1$ or for analytic sets).

\section{Sketch of the proof.} 

For Theorem \ref{Theorem1}, take $x_1,\dots,x_N$ a $(n,2\delta)$ separated set in $\mathrm{supp}(\nu)$, then $N$ is almost $\exp( n h_\nu(f))$.  For $x_i$, we show that the dynamical ball $B_n(x_i,\delta)$ contains $B(x_i, e^{- \chi_1 n})$ (shrinking that ball a little). Let $1_A$ denote the indicator function of $A$. Since the $B_n(x_i,\delta)$ are distinct and by definition of $d_1$, we have:
\begin{align*}
(d_1)^n& \simeq \int (f^n)_*(S )\wedge \omega \geq \sum_{i} \int (f^n)_*( 1_{B_n(x_i,\delta)} S )\wedge \omega\\
       & \geq \sum_{i} \int  1_{B_n(x_i,\delta)} S \wedge (f^n)^*(\omega) \geq \sum_{i} \int  1_{B(x_i, e^{- \chi_1 n})} S \wedge (f^n)^*(\omega) . 
\end{align*}
Assume that we can approximate the form $\omega$ by the current of integration along the unstable manifolds : say for example that locally the coordinates are $(x,y)$ so that $W^1_{(x,y)} = (x,y)+ t(1,0)$ and $W^2_{(x,y)} = (x,y)+ t(0,1)$ for all $(x,y)$ ($W^i_{(x,y)}$ is the unstable manifold associated to the Lyapunov exponent $\chi_i$). Then $\omega = idx \wedge d\bar{x} + idy \wedge d\bar{y}$ and $(f^n)^*(\omega) \simeq e^{2 \chi_1 n}idx \wedge d\bar{x}+ e^{2 \chi_2 n} idy \wedge d\bar{y} \geq  e^{2 \chi_2 n} \omega.$
So we have, using that $\int_{B(x,r)} S \wedge \omega \geq r^{\bar{d}_S}$:
$$(d_1)^n \geq \sum_{i} e^{2 \chi_2 n} \int  1_{B(x_i, e^{- \chi_1 n})} S \wedge  \omega \geq N e^{2 \chi_2 n} e^{- \bar{d}_S \chi_1 n}  . $$
The result then follows by taking the logarithm and letting $n \to \infty$.\\

For Theorem \ref{Theorem2}, we consider again the points $x_i$ which are $(n,2\delta)$ separated. By Brin-Katok's formula, we know that $\nu(B_n(x,\delta))\simeq \exp(-n h_\nu(f))$. We take points $y_1,\dots,y_L$ in  $B_n(x_i,\delta)$ inductively: we take $y_k$ in $B_n(x_i,\delta) \backslash \cup_{j<k} B(y_j, 2\exp(-\chi_1 n))$. As $\nu(B(y_j, 2\exp(-\chi_1 n))) \leq 2^{\underline{d}_\nu } e^{-\underline{d}_\nu \chi_1 n}\simeq  e^{-\underline{d}_\nu \chi_1 n} $, we deduce from a measure argument that we can take $L \geq \exp((\underline{d}_\nu\chi_1-h_\nu(f))n)$ and the points $y_1,\dots,y_L$ are $2\exp(-\chi_1 n)$-separated. 
Then we apply the arguments of the proof of Theorem \ref{Theorem1} to all the $y_j$ (for all $x_i$):
$$d_1^n\geq \sum_{i,j} \int (f^n)_*(1_{B(y_j,\exp(-\chi_1 n))} S)\wedge \omega \geq N.L e^{2\chi_2 n } e^{-\chi_1 n \bar{d}_S}. $$  
Again, the result follows by taking the logarithm and letting $n \to \infty$. \\

In order to make the previous arguments correct, we need first to use Pesin's theory to find sets where the estimates are uniform. But the main difficulty lies in the fact that the size of the unstable manifolds is not defined uniformly and they do not vary as a smooth foliation. So we take two smooth transverse foliations by lines near a point and we show using a theorem of graph transform that their behaviors are controlled under pull back by $f^n$. Finally, we have to distinguish the cases $\chi_1>\chi_2$ and $\chi_1=\chi_2$ (in the second case, graph transforms are not controlled but not necessary). Finally, the meromorphicity is not an issue as the balls we will consider do not meet the indeterminacy set. We will just need an approximation's argument to get the term in $d_1$ in the spirit of \cite{DS4}.

\section{Oseledec's Theorem and Pesin's theory} 
 We recall some facts on Pesin's theory and we define measurable sets where the estimates are uniform. The results of this section are taken from \cite{DT1}[p94--100]. Let $\nu$ be an ergodic, $f$-invariant measure with $\log d(x, I) \in L^1(\nu)$ such that $\chi_1 \geq \chi_2 >0$. Let $\widehat{X}$ be the natural extension of $X$:
 $$\widehat{X}:= \left\{ \hat{x}=(\dots, x_{-n},\dots, x_0, \dots, x_n, \dots )\in (X)^\Zz, \ \forall n, \ f(x_n)=x_{n+1} \right\}.$$
 In that space, $f$ induces a map $\hat{f}$ given by the left-side shift. Let $\Pi:\widehat{X}\to X$ be the canonical projection $\Pi(\hat{x})=x_0$ (so $f \circ \Pi  = \Pi \circ \hat{f}$). Then there exists a unique invariant probability measure $\hat{\nu}$ on $\widehat{X}$ such that $\Pi_* \hat{\nu}=\nu$. In that space, we only consider orbits that do not meet $\Acal$:
 $$\widehat{X}^*:=\left\{ \hat{x} \in \widehat{X}, \ \forall n, \ x_n\notin \Acal \right\}.  $$ 
Since $\chi_1\geq \chi_2 >0$, we have that $\log |\mathrm{Jac}(f)|\in L^1(\nu)$. In particular, $\widehat{X}^*$ is a set of full measure and we are in the settings of \cite{DT1}[p94--100]. We shall use the results of that paper keeping the same notations: $\tau_x$, $\varepsilon_0$, $f_x$, $f^n_x$, $f^{-n}_{\hat{x}}$, $D\hat{f}{x}$ and replace $\mu$ by $\nu$.  Then we have Oseledec's Theorem and Pesin $\gamma$-reduction's Theorem (we will choose a very small $\gamma$ later). Again, we keep the notations of \cite{DT1} so $q=1$ or $2$ depending if $\chi_1=\chi_2$ or not and $C_\gamma$ is the tempered map that gives the change of basis to put the differential in its simplified form. We now write $g_x$ the expression of $f_x$ in the charts given by $C_\gamma$:
$$g_x:= C_\gamma^{-1}(\hat{f}(\hat{x}))\circ f_x \circ C_\gamma(\hat{x})  $$
where $\Pi(\hat{x})=x$ (in \cite{DT1}, the notation used is $g_{\hat{x}}$ instead of $g_x$).  

In particular, \cite{DT1}[Proposition 8, Proposition 9] apply in our setting.
 Observe that, for Proposition 9, we denote 
\begin{align*}
g^{-1}_{\hat{x}}&:= C_\gamma^{-1}(\hat{f}^{-1}(\hat{x}))\circ f^{-1}_{\hat{x}} \circ C_\gamma(\hat{x})\\
                &=  C_\gamma^{-1}(\hat{f}^{-1}(\hat{x}))\circ f^{-1}_{x_{-1}} \circ C_\gamma(\hat{x}),
\end{align*}
in particular $g^{-1}_{\hat{x}} \neq (g_{x})^{-1}$ (this is the reason we change the notations in order to avoid confusion). Lemma 10 of \cite{DT1} also applies in our case, though in order to simplify the computations, we take here $d(x_n,\Acal)\geq V(\hat{x})e^{-\frac{\gamma|n|}{p}}$ where $p$ is the integer  given by \cite{DT1}[Proposition 9].

 Finally, let $\widehat{Y}$ be the set of points where the conclusions of the previous theorems, propositions and lemma are checked ($\hat{\nu}(\widehat{Y})=1$). \\

We now define sets where estimates are uniform. First, let:
$$ \Lambda_{k_0}:= \left\{x\in \Lambda, \ \forall n \geq k_0, \ \frac{\log(S\wedge \omega(B(x,e^{-\chi_1n- 15 \gamma n})))}{\log e^{-\chi_1n- 15 \gamma n}} \leq \bar{d}_S +\gamma \right \}.$$
Since $\bigcup_{k_0 \geq 0}\Lambda_{k_0}= \Lambda$ by (\ref{hypotheseS}) and $\nu(\Lambda)=1$, we choose $k_0$ large enough such that $\nu(\Lambda_{k_0}) \geq 3/4$. 

Consider now:
$$X_{\delta, n}:= \left\{ x, \ \nu(B_n(x,\delta)) \leq e^{-h_{\nu}(f)n +\gamma n}    \right\}. $$
Brin-Katok's formula gives that for $\delta$ small enough:
\begin{align*}
\frac{4}{5} &\leq \nu \left(\left\{x, \ \underline{\lim} -\frac{1}{n} \log \nu(B_n(x,\delta))\geq h_\nu(f) -\frac{\gamma}{2} \right\}\right) \\
            &\leq \nu \left(  \bigcup_{n_0} \bigcap_{n \geq n_0} X_{\delta,n}   \right).
\end{align*} 
In particular, we choose $n_0$ large enough so that:
$$\nu\left(\bigcap_{n\geq n_0} X_{\delta,n_0}\right) \geq \frac{3}{4}.$$
Let 
$$\widehat{Y}_{\alpha_0}:=\left\{ \hat{x}\in \widehat{Y}, \ \alpha_0 \leq \| C^{\pm1}_\gamma(\hat{x}) \|\leq \frac{1}{\alpha_0} \ \mathrm{and} \ V(\hat{x}) \geq \alpha_0 \right\}$$
where we have $d(x_n,\Acal)\geq V(\hat{x})e^{-\frac{\gamma|n|}{p}}$ for all $n\in \Zz$. Since $\hat{\nu}(\widehat{Y})=1$, we can choose $\alpha_0$ small enough so that $\hat{\nu}(\widehat{Y}_{\alpha_0})\geq 3/4$. 

Let $A:= \Pi(\widehat{Y}_{\alpha_0}) \cap \left(\bigcap_{n\geq n_0} X_{\delta,n_0}\right) \cap \Lambda_{k_0}$. Since $\Pi_*(\hat{\nu})=\nu$, we have that $\nu(\Pi(\widehat{Y}_{\alpha_0})) \geq \hat{\nu}(\widehat{Y}_{\alpha_0}) \geq 3/4$, hence $\nu(A)\geq 1/4$. \\

We will also need the graph transform Theorem of the first author (\cite{DT1}[Théorème p.100]). We state it here for the reader's convenience after redefining some notations. Consider the following map $g$ in the ball $B( 0,r)\subset \Cc^2$:
$$g(X,Y)=(AX + R(X,Y), BY + U(X,Y))$$
with $E_1:=\{(X,0)\}$, $E_2:=\{(0,Y)\}$  and $A$, $B$ complex numbers. Assume that $g(0,0)=(0,0)$ and consider $\|Dh(Z) \|:=\max (\| D
R(Z) \|, \| D U (Z) \|) $ in the ball $B(0,r)$. We also assume that   $\gamma \leq \| A \| \leq \| B \|(1 - \gamma)$ (observe that $\| B^{-1} \|^{-1} = \|B\|=|B|$ as $B\in \Cc$). Let $\{ ( \Phi(Y),Y), \ Y \in D \}$ be a 
graph in $B(0,r)$ for $Y\in D \subset E_2$ such that $\| \Phi(Y_1) - \Phi(Y_2) \| \leq \gamma_0 \| Y_1 - Y_2
\|$. The theorem gives conditions on $\|Dh \|$,
$\gamma$ and $\gamma_0$ so that the image of that graph by $g$ is again a
graph that satisfies the same control.
\begin{theorem}\label{p.100}
If  $\|Dh \|_{\infty} \| B^{-1} \| (1 + \gamma_0) <1$ then the image by $g$ of the previous graph is a graph above $\pi_0(g(\mbox{graph of \ }
\Phi))$ where $\pi_0$ is the projection on $E_2$. Furtheremore, if we denote $(\Psi(Y),Y)$ that new graph, then:
$$\| \Psi(Y_1) - \Psi(Y_2) \| \leq \frac{\| A \| \gamma_0 + \|Dh \|(1+
  \gamma_0)}{\| B^{-1} \|^{-1} -\|Dh \|(1+ \gamma_0)}  \| Y_1 - Y_2
  \|$$
which is less than $\gamma_0 \| Y_1 - Y_2 \|$ if $\|Dh \| \leq
  \epsilon(\gamma_0, \gamma)$.

Finally, if furthermore  $B(0 , \alpha) \subset D$ and $\| \Phi(0) \|
\leq \beta$, then $\pi_0(g(\mbox{graph of \  } \Phi))$ contains  $B(0,
(\| B^{-1} \|^{-1} - \|Dh \| (1 + \gamma_0)) \alpha - \| A \| \beta -
\|Dh \| \beta - \|
D^2 g \|_{B(0,r)} \beta^2)$ and $\| \Psi(0) \| \leq   (1+ \gamma_0)( \| A \| \beta + \|Dh \|  \beta +
\|D^2g \|_{B(0,r)} \beta^2) $ (if $\|Dh \| \leq \epsilon(\gamma_0, \gamma)$).
\end{theorem}

\section{Proof of Theorem \ref{Theorem1}}
In what follows, we take $n \geq \max(n_0, k_0)$. If $x\in A$, we have that $\nu(B_n(x,\delta))\leq e^{-h_{\nu}(f)n+\gamma n}$. Hence we can find a $(n,\delta)$ separated set $(x_1,\dots, x_N)$ in $A$ with 
$$N \geq \frac{1}{4} e^{h_{\nu}(f)n-\gamma n}.$$
 The Bowen balls $B_n(x_i,\delta/2)$ are distinct. In what follows, we fix $x$ as one of the $x_i$ and we give some properties around $x$ for 
$g_x^n:=g_{f^{n-1}(x)}\circ \dots \circ g_x$.
\begin{lemma}\label{lemma3}
For $l = 0,\dots,n$, we have:
$$g_x^l(B(0,e^{-\chi_1 n -12 \gamma n})) \subset B(0, e^{-\chi_1 (n-l)-12 \gamma n+ 2 \gamma l}),$$
in particular,
\begin{align*}
 g_x^n(B(0,e^{-\chi_1 n -12 \gamma n})) &\subset B(0, e^{-10 \gamma n}) \\
B(x,e^{-\chi_1 n-13\gamma n})&\subset B_n(x, \delta/2). 
\end{align*}
\end{lemma} 
\noindent \emph{Proof.} We prove the property by induction on $l$, the case $l=0$ is clear. So, assume it is true for $l \leq n-1$. Then for $l+1$:
\begin{align*}
g_x^{l+1}(B(0,e^{-\chi_1 n -12 \gamma n}))&= g_{f^{l}(x)}( g_x^l(B(0,e^{-\chi_1 n -12 \gamma n}))) \\
                                           &\subset  g_{f^{l}(x)}( B(0, e^{-\chi_1 (n-l)-12 \gamma n+ 2 \gamma l})). 
\end{align*}
Let $y \in g_{f^{l}(x)}( B(0, e^{-\chi_1 (n-l)-12 \gamma n+ 2 \gamma l}))$, then we have $y = g_{f^{l}(x)}(z) $ for some $z \in B(0, e^{-\chi_1 (n-l)-12 \gamma n+ 2 \gamma l})$. Then:
\begin{equation}\label{lemme3eq1}
 \|g_{f^{l}(x)}(z)-g_{f^{l}(x)}(0) \| \leq  \| Dg_{f^{l}(x)} \|_{ B(0, e^{-\chi_1 (n-l)-12 \gamma n+ 2 \gamma l})} \|z\|.
\end{equation}
Using \cite{DT1}[Proposition 8] gives that for:
$$\| w\|\leq \frac{\varepsilon_0 d(f^l(x),\Acal)^p}{\|C_\gamma(\hat{f}^l(\hat{x}))\|}$$
(recall that $\varepsilon_0$ is the radius of the ball $B(0, \varepsilon_0)$ on which the maps $\tau_x$ are defined), we have:
\begin{align*}
 \| Dg_{f^{l}(x)}(w) - Dg_{f^{l}(x)}(0)\| & = \| Dh(w) \| \\
                                          & \leq \tau \| C_\gamma(\hat{f}^{l+1}(\hat{x}))^{-1}\|\times  \| C_\gamma(\hat{f}^{l}(\hat{x}))\|^2 d(f^l(x),\Acal)^{-p} \|w\|\\
																					&\leq \tau e^{(l+1)\gamma} e^{2l \gamma}\frac{1}{\alpha_0^3} \frac{1}{\alpha_0^p}e^{\frac{\gamma n p }{p}} \|w\| 
\end{align*}
where we used that $\|C_\gamma^\pm\| $ is tempered (see \cite{KH}[Lemma S.2.12, p. 668]) and $\hat{x}\in \widehat{Y}_{\alpha_0}$. Now, we see that:
\begin{align*}
\frac{\varepsilon_0 d(f^l(x),\Acal)^p}{\|C_\gamma(\hat{f}^l(\hat{x}))\|} &\geq \varepsilon_0 e^{-\gamma l}\alpha_0^2 e^{-\gamma l} \\
                                                                   &\geq e^{-\chi_1(n-l)-12 \gamma n +2 \gamma l}.
\end{align*}
Hence, we can apply the previous bound of $ \| Dg_{f^{l}(x)}(w) - Dg_{f^{l}(x)}(0)\| $ to $\|w\| \leq e^{-\chi_1(n-l)-12 \gamma n +2 \gamma l}$:
\begin{align*}
 \| Dg_{f^{l}(x)}(w) - Dg_{f^{l}(x)}(0)\| \leq \frac{C}{\alpha_0^{p+3}}e^{-6\gamma n} 
\end{align*}
Using \cite{DT1}[Théorème de $\gamma$-r\'eduction de Pesin, Proposition 8] gives that:
$$\|Dg_{f^l(x)}(0)\|\leq e^{\chi_1 +\gamma}. $$
In particular:
$$\|Dg_{f^{l}(x)}(w)\| \leq e^{\chi_1 +\gamma}+ \frac{C}{\alpha_0^{p+3}}e^{-6\gamma n} \leq e^{\chi_1 +2 \gamma} $$ 
for $n$ large. Then (\ref{lemme3eq1}) becomes:
\begin{align*}
 \|g_{f^{l}(x)}(z)-g_{f^{l}(x)}(0) \| &\leq  e^{\chi_1 + 2 \gamma}\|z\| \\
                                      &\leq e^{\chi_1 + 2 \gamma}e^{-\chi_1(n-l)-12\gamma n+2\gamma l}\\
                                      &\leq e^{-\chi_1(n-l-1)-12\gamma n +2 \gamma (l+1)},
\end{align*}
and the induction is proved.

Now, as $g_x^l:=g_{f^{l-1}(x)}\circ \dots \circ g_x$, we can write, using the notations $f_{x}$ and $f^l_{x}$ of \cite{DT1}:
\begin{align*}
g_x^l&= C^{-1}_\gamma(\hat{f}^l(x))\circ f_{f^{l-1}(x)} \circ \dots \circ f_x \circ C_\gamma(\hat{x})\\
     &= C^{-1}_\gamma(\hat{f}^l(x))\circ \tau_{f^{l}(x)}^{-1} \circ f^l \circ \tau_{x}\circ C_\gamma(\hat{x}).
\end{align*}
We can assume that the first derivatives of $\tau_x^{-1}$ are uniformly bounded (independently of $x$) and we use that $\|C_\gamma(\hat{x})^{-1}\|\leq \alpha_{0}^{-1}$, hence:
$$ \left(\tau_{x}\circ C_\gamma(\hat{x})\right)^{-1}B(x,e^{-\chi_1 n -13 \gamma n})\subset B(0,e^{-\chi_1 n -12 \gamma n}). $$
Hence:
\begin{align*}
 f^l\left(B(x,e^{-\chi_1 n -13 \gamma n})\right) & \subset f^l \circ \tau_x \circ C_{\gamma}(\hat{x})\left(B(0,e^{-\chi_1 n -12 \gamma n})  \right)  \\
                                                 & \subset \tau_{f^{l}(x)}\circ C_\gamma(\hat{f}^l(x)) \circ g_x^l(B(0,e^{-\chi_1 n -12 \gamma n}))\\
                                                 &\subset  \tau_{f^{l}(x)}\circ C_\gamma(\hat{f}^l(x))(B(0,e^{-\chi_1 (n-l) -12 \gamma n+2\gamma l}))
 \end{align*}
where the last inclusion follows from the previous induction. Hence, using that $C_\gamma$ is tempered and the bound on the derivatives of $\tau_x$ by a constant $C$ implies:
$$  f^l\left(B(x,e^{-\chi_1 n -13 \gamma n})\right) \subset B(f^l(x), \frac{Ce^{\gamma l }}{\alpha_0}e^{-\chi_1(n-l)-12\gamma n+2 \gamma l }) \subset B(f^l(x),\frac{\delta}{2}),$$
for $l=0,\dots,n$ (we take $n$ large so that, for example, $2\exp(-\gamma n ) \leq \delta$). The lemma follows. \hfill $\Box$ \hfill\\

\subsection{Case where $\chi_1 > \chi_2$}
We consider now the case where $\chi_1 > \chi_2$. Recall that $x$ denotes one the $x_i$. In particular, 
$$S\wedge \omega (B(x,e^{-\chi_1 n -15 \gamma n}))\geq e^{(-\chi_1 n-15 \gamma n)(\bar{d}_S+\gamma)}. $$
We first pull that inequality back to $\Cc^2$.
Let $\gamma_0$ be a constant such that 
\begin{align}\label{gamma0}
e^{\chi_2+\gamma} + \gamma_0 e^{\chi_1+\gamma} &<  e^{\chi_2 + 2 \gamma}  \nonumber  \\
\frac{e^{\chi_2-\gamma} - \gamma_0 e^{\chi_1+\gamma}}{\sqrt{1+\gamma_0^2}}&>  e^{\chi_2 - 2 \gamma},
\end{align}
 (see the end of the proof of Lemma \ref{lemma9}) and $(1+\gamma_0)\leq e^{\chi_1-\chi_2 - 2 \gamma}$ (see the end of the proof of Lemma \ref{Forward graph transform}). Let $D_1$ and $D_2$ be two lines in $\Cc^2$ passing through $0$ and making an angle greater than $\gamma_0$. Assume furthermore that $D_1$ and $D_2$ make an angle with $C_\gamma^{-1}(E_{1}(\hat{x}))$ smaller than $\gamma_0$ ($E_1(\hat{x})$ is defined in \cite{DT1}[Th\'eor\`eme d'Osedelec]). Consider the positive $(1,1)$-form $\beta_1$ on $\Cc^2$ obtained by integrating over all the currents of integration on lines parallel to $D_1$:
$$\forall \psi \in \Dcal^{1,1}, \ \int \beta_1\wedge \psi = \int_{D_1^{\bot}} \int_{D_1(z)} \psi idz \wedge d\bar{z} $$
where $D_1(z)$ is given by the translation of $D_1$ passing through $z\in D_1^{\bot}$. Similarly, we define $\beta_2$. Let $\beta$ denote the positive smooth $(1,1)$ form $ \beta_1+\beta_2$. We have in a ball $B(x, C \varepsilon_0) \supset B(x,e^{- \gamma n})$:
$$ (\tau_x \circ C_\gamma(\hat{x}))_* \beta \geq C(\alpha_0, \gamma_0) \omega $$
for some constant $C(\alpha_0, \gamma_0)$. We denote in what follows $\Scal := (\tau_x \circ C_\gamma(\hat{x}))^*(S)$. Then:
\begin{align*}
&(\tau_x \circ C_\gamma(\hat{x}))^*(S\wedge \omega)\left((\tau_x \circ C_\gamma(\hat{x}))^{-1}\left(B\left(x,e^{-\chi_1 n -15 \gamma n}\right)\right)\right)=\\
&S \wedge \omega \left(B\left(x,e^{-\chi_1 n -15 \gamma n}\right)\right)\geq e^{(-\chi_1 n -15 \gamma n)(\bar{d}_S+\gamma)}.
\end{align*}
So:
$$\Scal \wedge \beta \left((\tau_x \circ C_\gamma(\hat{x}))^{-1}\left(B\left(x,e^{-\chi_1 n -15 \gamma n}\right)\right)\right) \geq C(\alpha_0, \gamma_0) e^{(-\chi_1 n -15 \gamma n)(\bar{d}_S+\gamma)}.$$
So we can assume that:
$$\Scal \wedge \beta_1 \left((\tau_x \circ C_\gamma(\hat{x}))^{-1}\left(B\left(x,e^{-\chi_1 n -15 \gamma n}\right)\right)\right) \geq \frac{C(\alpha_0, \gamma_0)}{2} e^{(-\chi_1 n -15 \gamma n)(\bar{d}_S+\gamma)}.$$
Finally, we can assume that the first order derivatives of the $(\tau_x)^{-1}$ are uniformly bounded so that:
$$ (\tau_x \circ C_\gamma(\hat{x}))^{-1}\left(B\left(x,e^{-\chi_1 n -15 \gamma n}\right)\right)\subset B\left(0,e^{-\chi_1 n -14 \gamma n}\right).$$
We deduce:
\begin{equation}
\label{eqchi_1}
\Scal \wedge \beta_1(B\left(0,e^{-\chi_1 n -14 \gamma n}\right) ) \geq  \frac{C(\alpha_0, \gamma_0)}{2} e^{(-\chi_1 n-15 \gamma n)(\bar{d}_S+\gamma)}.
\end{equation}

We now do some graph transforms with the leaves $D_1(z)$ composing $\beta_1$. We start with some leaf that intersects $B(0,e^{-\chi_1 n-12 \gamma n})$. That leaf is a graph over $C^{-1}_{\gamma}E_1(\hat{x})$. We denote that graph by $(X,\Phi_0(X))$ for $|X|<  e^{-\chi_1n-7\gamma n}$. As the angle between $D_1$ and $C^{-1}_{\gamma}E_1(\hat{x})$ is $\leq \gamma_0$ then $ \|\Phi_0(X)-\Phi_0(X') \|\leq \gamma_0 \|X-X' \|$. Furthermore, since the graph intersects $B(0,e^{-\chi_1 n-12 \gamma n})$, we have that $\|\Phi_0(0)\| \leq e^{- \chi_1 n -10 \gamma n}$ ($\gamma_0 \ll n$). We will consider the images of that graph by $g^l_{x}$ for $1\leq l\leq n-1$ and show that at each step we still have a graph with a good control. 
\begin{lemma}[Forward graph transform]\label{Forward graph transform}

For $0\leq l\leq n-1$,  let $(X, \Phi_l(X))$ be a graph over $C^{-1}_\gamma(\hat{f}^l(\hat{x}))E_1(\hat{f}^l(\hat{x}))$ for $|X|< e^{-\chi_1(n-l)-7\gamma n -2 \gamma l}$. Assume that 
$$\mathrm{Lip}(\Phi_l) \leq \gamma_0 \  \mathrm{and} \ |\Phi_l(0)| \leq e^{-\chi_1(n-l)-10 \gamma n}.$$ 
Then the image of that graph by $g_{f^l(x)}$ is a graph over $C^{-1}_\gamma(\hat{f}^{l+1}(\hat{x}))E_1(\hat{f}^{l+1}(\hat{x}))$ at least for $|X|< e^{-\chi_1(n-l)-7\gamma n -2 \gamma (l+1)}$ and
$$\mathrm{Lip}(\Phi_{l+1}) \leq \gamma_0 \ \mathrm{and} \ | \Phi_{l+1}(0)| \leq e^{-\chi_1(n-l-1)-10 \gamma n}$$
\end{lemma}
\noindent \emph{Proof.}  Pesin's Theorem \cite{DT1} gives that in 
$$\Cc \oplus \Cc= C^{-1}_\gamma(\hat{f}^l(\hat{x}))E_1(\hat{f}^l(\hat{x})) \oplus C^{-1}_\gamma(\hat{f}^l(\hat{x}))E_2(\hat{f}^l(\hat{x})) ,$$
 the map $g_{f^l(x)}$ is given by:
$$g_{f^l(x)}(X,Y)= (AX + R(X,Y), B Y +U(X,Y))  $$
 with $A= A^1_\gamma(\hat{f}^l(\hat{x}))$, $ B= A^2_\gamma(\hat{f}^l(\hat{x}))$ satisfying:
$$\gamma \leq e^{\chi_2-\gamma}\leq \|B\| \leq e^{\chi_1-\gamma}(1- \gamma) \leq (1-\gamma) \|A^{-1}\|^{-1} $$
since $\gamma$ is small with respect to $\chi_1$ and $\chi_2$. We want to apply Theorem \ref{p.100} (observe that we exchange here the role of $X$ and $Y$). In order to do that, we need to control $\|D h(z)\|=\max (\|D R(z)\| ,\| D U(z) \|)$. 
Observe that for $(X,\Phi_l(X))$ with $|X|\leq e^{-\chi_1(n-l)-7\gamma n -2 \gamma l}$ then:
\begin{align*}
\|(X,\Phi_l(X)) \| &\leq |X| + |\Phi_l(0)| + |\Phi_l(0)-\Phi_l(X) | \\
     &\leq  e^{-\chi_1(n-l)-7\gamma n -2 \gamma l}+e^{-\chi_1(n-l)-10 \gamma n} + \gamma_0 |X|\\
     &\leq 3 e^{-\chi_1 (n-l) -7 \gamma n -2 \gamma l}\leq e^{-5\gamma n} = r
\end{align*}
with the notations of Theorem \ref{p.100}. Hence, by \cite{DT1}[Proposition 8], we have:
$$\|Dh(z)\|\leq C \|C_\gamma(\hat{f}^{l+1}(\hat{x}))^{-1} \| \|C_\gamma(\hat{f}^{l}(\hat{x})) \|^2 d(f^l(x),\Acal)^{-p}\|z\| $$
for $\|z\|\leq \varepsilon_0 \frac{d(f^l(x),\Acal)^p}{\|C_\gamma(\hat{f}^{l}(\hat{x})) \|}$. As in the proof of Lemma \ref{lemma3}, we check that:
\begin{equation*}
 \varepsilon_0 \frac{d\left(f^l(x),\Acal \right)^p}{\|C_\gamma(\hat{f}^{l}(\hat{x})) \|}\geq \varepsilon_0 e^{-\gamma l}\alpha_0 \alpha_0 e^{-\gamma l} \geq \varepsilon_0 \alpha_0^2e^{-2\gamma n} \geq e^{-3 \gamma n}. 
\end{equation*}
In particular, we can apply the previous estimate to $\|z\| \leq e^{-5 \gamma n }$:
\begin{align}\label{estimeDh}
\|Dh(z)\|&\leq \frac{C}{\alpha_0^3} e^{\gamma (l+1)+ 2\gamma l}   \frac{1}{\alpha_0^p}e^{\gamma l} \|z\|  \nonumber\\
         & \leq \frac{C}{\alpha_0^{p+3}} e^{4\gamma n} \|z\| \leq \frac{C}{\alpha_0^{p+3}} e^{4\gamma n} e^{-5\gamma n} \leq \frac{C}{\alpha_0^{p+3}}e^{-\gamma n}
\end{align}
which is arbitrarily small for $n$ large enough (hence $\leq \varepsilon(\gamma_0,\gamma)$). In particular $\|D h\|$ is sufficiently small to provide:
$$\|D h\|\|A^{-1}\|(1+\gamma_0) <1. $$
Let $\alpha:=e^{-\chi_1(n-l)-7\gamma n -2 \gamma l}$ and $\beta:= e^{-\chi_1(n-l)-10 \gamma n}$, then we can apply Theorem \ref{p.100} as $(X,\Phi_l(X))$ is a graph over $B(0,\alpha)$ and $|\Phi_l(0)| \leq \beta$. We deduce that $(X, \Phi_{l+1}(X))$ is a graph at least for 
\[|X| < (\|A^{-1}\|^{-1} - \|D h\| (1+\gamma_0))\alpha - \|B\|\beta - \|D h\|\beta -  \| D^2g_{f^l(x)}\|_{B(0,r)} \beta^2 \]
where $r=e^{-5 \gamma n}$. So it is a graph for:
\begin{align*}
|X|\leq & \left(e^{\chi_1-\gamma} - \frac{C}{\alpha_0^{p+3}}e^{-\gamma n}(1+\gamma_0) \right)e^{-\chi_1(n-l)-7\gamma n -2 \gamma l} -  \\ 
&  e^{-\chi_1(n-l)-10 \gamma n}\left(e^{\chi_2 +\delta}+ \frac{C}{\alpha_0^{p+3}}e^{-\gamma n}+ \| D^2g_{f^l(x)}\|_{B(0,r)}e^{-\chi_1(n-l)-10 \gamma n} \right). 
\end{align*}
Using \cite{DiDu}[Lemma 2], as in \cite{DT1}[Proposition 9], implies that:
$$ \|D^2g_{f^l(x)}\|_{B(0,r)} \leq e^{5 \gamma n} $$ 
 since the $C^{\pm 1}_{\gamma}$ are tempered. Hence, we have a graph for:
  \begin{align*}
|X| \leq & \Big(  e^{\chi_1-\gamma} - \frac{C}{\alpha_0^{p+3}}e^{-\gamma n}(1+\gamma_0) -  \\ 
&  e^{-3 \gamma n +2\gamma l}( e^{\chi_2 +\delta}- \frac{C}{\alpha_0^{p+3}}e^{-\gamma n}-e^{-\chi_1(n-l)-5\gamma n}) \Big)   e^{-\chi_1(n-l)-7\gamma n -2 \gamma l}.\end{align*}
Hence we have a graph at least for:
  \begin{align*}
|X| &\leq e^{\chi_1 -2\gamma} e^{-\chi_1(n-l)-7\gamma n -2\gamma l}\\
    &=e^{-\chi_1(n-l-1)-7 \gamma n -2\gamma (l+1)} 
\end{align*}
 which is what we want. 
 
 To conclude, we have to bound $|\Phi_{l+1}(0)|$, for that we still apply  Theorem \ref{p.100}:
 \begin{align*}
 |\Phi_{l+1}(0)|&\leq (1+\gamma_0)(\|B\|+ \|D h\| +\|D^2g_{f^l(x)}\|_{B(0,r)} \beta)\beta \\
                  &\leq (1+\gamma_0)(e^{\chi_2+\gamma}+ \frac{C}{\alpha_0^{p+3}}e^{-\gamma n} + e^{5\gamma n} e^{-\chi_1(n-l)-10\gamma n})\beta\\
                  &\leq (1+\gamma_0)e^{\chi_2+2\gamma}\beta\\
                  &\leq e^{\chi_1}\beta=e^{-\chi_1(n-l-1)-10\gamma n} 
 \end{align*}  
 where the last inequality follows from $(1+\gamma_0)\leq e^{\chi_1-\chi_2 - 2 \gamma}$. \hfill $\Box$ \hfill \\

Observe that we have the following corollary of the proof (see (\ref{estimeDh})) that we will use later, it does not use the fact that the Lyapunov exponents are distinct:
\begin{corollary}\label{estimedh}
With the notations of the proof of Lemma \ref{Forward graph transform}, we have that for $\|z \| \leq e^{-3\gamma n}$:
$$ \|Dh(z)\| \leq \frac{C}{\alpha_0^{p+3}}e^{4 \gamma n } \|z\|. $$
\end{corollary}
Now, we apply inductively the lemma for $l=0\dots n-1$ and we only keep at each step the part over $B(0,e^{-\chi_1(n-l)-7\gamma n -2 \gamma l})$. We get that $(X,\Phi_n(X))$ is a graph for $|X|\leq e^{-9\gamma n}$ with $\mathrm{Lip}(\Phi_n)\leq \gamma_0$ and $|\Phi_n(0)|\leq e^{-10\gamma n}$. Observe that, at each step, we drop no part of $g^n_x(D_1(z) \cap B(0,e^{-\chi_1 n -12\gamma n}))$. Indeed, Lemma  \ref{lemma3} gives 
 $$g^n_x(B(0,e^{-\chi_1 n-12\gamma n })) \subset B(0,e^{-\chi_1(n-l)-12\gamma n + 2 \gamma l}) $$
 and that is contained in the part over $|X|\leq e^{-\chi_1(n-l)-7\gamma n -2 \gamma l}$ that we keep. We now consider backward graph transforms. We pull back $B(0,e^{-7\gamma n})$ which is contained in $C_\gamma^{-1}(\hat{f}^n(\hat{x})) E_2(\hat{f}^n(\hat{x}))$ by $g^n_x$. Observe that this disk is a graph $(\Psi_n(Y),Y)$ with $\Psi_n\equiv 0$. 
\begin{lemma}[Backward graph transform]\label{Backward graph transform}

For $0\leq l\leq n-1$,  let $(\Psi_{n-l}(Y), Y)$ be a graph over $C^{-1}_\gamma(\hat{f}^{n-l}(\hat{x}))E_2(\hat{f}^{n-l}(\hat{x}))$ for $|Y|< e^{-\chi_2l-7\gamma n -2 \gamma l}$. Assume that 
$$\mathrm{Lip}(\Psi_{n-l}) \leq \gamma_0 \  \mathrm{and} \ \Psi_{n-l}(0) =0.$$ 
Then the image of that graph by $g^{-1}_{\hat{f}^{n-l}(\hat{x})}$ is a graph over $C^{-1}_\gamma(\hat{f}^{n-l-1}(\hat{x}))E_2(\hat{f}^{n-l-1}(\hat{x}))$ at least for $|Y|< e^{-\chi_2(l+1)-7\gamma n -2 \gamma (l+1)}$ and
$$\mathrm{Lip}(\Psi_{n-l-1}) \leq \gamma_0 \ \mathrm{and} \ \Psi_{n-l-1}(0)=0.$$
\end{lemma}  
 \noindent \emph{Proof.} We apply \cite{DT1}[Proposition 9]. First, we know that $g^{-1}_{\hat{f}^{n-l}(\hat{x})}(w)$ is well defined for $\|w\| \leq \varepsilon_0 d(f^{n-l-1}(x), \Acal)^{p}/\| C_\gamma(\hat{f}^{n-l}(\hat{x}))\|$. Using again that $C_\gamma$ is tempered and that $d(f^{n-l-1}(x), \Acal)\geq \alpha_0 e^{-\gamma (n-l-1 ) /p }$ imply:
 $$\frac{\varepsilon_0 d(f^{n-l-1}(x), \Acal)^{p}}{\| C_\gamma(\hat{f}^{n-l}(\hat{x}))\|} \geq \varepsilon_0 e^{-\gamma (n-l+1)} \alpha_0^{p+1} e^{-\gamma n} \geq \varepsilon_0 e^{-2\gamma n} \alpha_0^{p+1}. $$
 In particular $g^{-1}_{\hat{f}^{n-l}(\hat{x})}(w)$ is well defined on $B(0,e^{-6\gamma n})$ which contains the graph $(\Psi_{n-l}(Y),Y) $. 
 
 We now do the graph transform. In
 $$\Cc \oplus \Cc= C^{-1}_\gamma(\hat{f}^{n-l}(\hat{x}))E_1(\hat{f}^{n-l}(\hat{x})) \oplus C^{-1}_\gamma(\hat{f}^{n-l}(\hat{x}))E_2(\hat{f}^{n-l}(\hat{x})) ,$$
 the map $g^{-1}_{\hat{f}^{n-l}(\hat{x})}$ is given by:
$$g^{-1}_{\hat{f}^{n-l}(\hat{x})}(X,Y)= (AX + R(X,Y), B Y +U(X,Y))  $$
with 
\begin{align*}
A&= A^1_\gamma(\hat{f}^{n-l-1}(\hat{x}))^{-1} \\
B&= A^2_\gamma(\hat{f}^{n-l-1}(\hat{x}))^{-1}.
 \end{align*}
 In particular, we have indeed:
 $$\gamma \leq e^{-\chi_1-\gamma}\leq \|A\| \leq e^{-\chi_2-\gamma}(1-\gamma)\leq (1-\gamma)\|B^{-1}\|^{-1} $$
 since $\gamma$ is small compared to the exponents $\chi_i$. Now we  bound $\max(\|DR\|,\|DU\|)=\|Dh\|$ on $B(0,e^{-6 \gamma n})$:
 \begin{align}\label{estimeDh2}
 \|Dh(w)\|&\leq \tau \|C^{-1}_\gamma( \hat{f}^{n-l-1}(\hat{x}))   \| \|C_\gamma( \hat{f}^{n-l}(\hat{x}))   \|^2 d(f^{n-l-1}(x), \Acal)^{-p}\|w\|\\
          &\leq \tau \frac{1}{\alpha_0^3} e^{3\gamma n} \frac{1}{\alpha_0^p} e^{\gamma n} \|w\| \nonumber \\
          &\leq \tau \frac{1}{\alpha_0^{p+3}} e^{-2\gamma n}.  \nonumber
 \end{align}
We apply the results of Theorem \ref{p.100}.  In particular, that term is small enough to provide $\|Dh(w)\| \|B^{-1}\|(1+\gamma_0) <1$ so the image  of $(\Psi_{n-l}(Y),Y)$ by $g^{-1}_{\hat{f}^{n-l}(\hat{x})}$ is a graph $(\Psi_{n-l-1}(Y),Y)$. As $ \|Dh(w)\|\leq \varepsilon(\gamma_0,\gamma)$ ($n$ is large), we have $\mathrm{Lip}(\Psi_{n-l-1}) \leq \gamma_0$ and $ \Psi_{n-l-1}(0) =0$ since  $g^{-1}_{\hat{f}^{n-l}(\hat{x})}(0)=0$. It remains to see that $(\Psi_{n-l-1}(Y),Y)$ is a graph for $|Y|< e^{-\chi_2(l+1)-7\gamma n -2 \gamma (l+1)}$. Here $\alpha=e^{-\chi_2l-7\gamma n -2 \gamma l}$ and $\beta=0$ hence $(\Psi_{n-l-1}(Y),Y)$ is a graph for
$$|Y|<\alpha(\|B^{-1}\|^{-1}-\|Dh(w)\|(1+\gamma_0)). $$
 Then:
 \begin{align*}
 \alpha(\|B^{-1}\|^{-1}-\|Dh(w)\|(1+\gamma_0))&\geq  ( e^{-\chi_2-\gamma}-\tau \frac{1}{\alpha_0^{p+3}} e^{-2\gamma n}(1+\gamma_0))\alpha\\
                                              &\geq e^{-\chi_2-2\gamma}e^{-\chi_2l-7\gamma n -2 \gamma l}\\
                                              &\geq e^{-\chi_2(l+1)-7\gamma n -2 \gamma (l+1)},
                                               \end{align*}
 which concludes the proof.\hfill $\Box$ \\
 
Using that $\tau \alpha_0^{-(p+3)}\leq \exp(\gamma n)$, observe that we have the following corollary of the proof (see (\ref{estimeDh2})) that we will use later, it does not use the fact that the Lyapunov exponents are distinct:
\begin{corollary}\label{estimedh2}
With the notations of the proof of Lemma \ref{Backward graph transform}, we have that for $\|w \| \leq e^{-6\gamma n}$:
$$ \|Dh(w)\| \leq e^{5 \gamma n } \|w\|. $$
\end{corollary}

 Using the lemma inductively for $l=0,\dots, n-1$, keeping only at each step the part above $ |Y| \leq e^{-\chi_2l-7\gamma n -2 \gamma l}$, gives a graph $(\Psi_0(Y),Y)$ for $|Y|\leq e^{-\chi_2n-9\gamma n}$ with 
 $\mathrm{Lip}(\Psi_{0}) \leq \gamma_0$ and $ \Psi_{0}(0) =0$. We denote that graph by $W$. Recall that the form $\beta_1$ was defined by
$$\int \beta_1\wedge \psi = \int_{D_1^{\bot}} \int_{D_1(z)} \psi .$$
We want to decompose $\beta_1$ along $W$ instead. In order to simplify the notations, we still denote by $(X,Y)$  the coordinates given by $(D_1,D_1^{\bot})$. For that we need the lemma:
 \begin{lemma}\label{lemma6}
The set $W$ is a graph $(\Phi(Y),Y)$ above $D_1^{\bot}$ for $|Y|<e^{-\chi_2 n-10 \gamma n}$. Its Lipschtitz' constant is less than a constant of $\gamma_0$ denoted by $C(\gamma_0)$.
\end{lemma}
\noindent \emph{Proof.} We make a change of coordinates. Let $(e_1,e_2)$ denote the initial orthonormal basis (i.e. $e_1 \in C_\gamma^{-1}(\hat{x})E_1(\hat{x})$, $e_2 \in  C_\gamma^{-1}(\hat{x})E_2(\hat{x})$). Let $\alpha= \cos (\mathrm{angle}(e_1,D_1))$ and  $\beta= \sin (\mathrm{angle}(e_1,D_1))$ where $\mathrm{angle}(e_1,D_1) \leq \gamma_0$ by hypothesis. Then  $|\alpha|\geq \cos(\gamma_0)$ and $|\beta|\leq |\sin(\gamma_0)|\leq \gamma_0$. Consider the matrix
$$P:=\begin{pmatrix}
	\alpha &  -\beta \\
	\beta      &  \alpha   
\end{pmatrix}. $$
Let $(v_1,v_2)$ be the orthonormed basis with $v_1 \in D_1$, $v_2 \in D_1^\bot$ given by $P^{-1}(v_1,v_2)=(e_1,e_2)$.
 We show that $P^{-1} \circ (\Psi_0(Y),Y)$ is a graph for $|Y|<e^{-\chi_2 n-10 \gamma n}$. We apply for that Theorem \ref{p.100}. We write :
 $$ P(X,Y)=(AX+R(X,Y), BY+U(X,Y))$$
 so that $A=B=\alpha$, $R(X,Y)=\beta Y$, $U(X,Y)=-\beta X$. We do not have the relation $\gamma \leq \|A \| \leq \|B\|(1-\gamma)$ but it is only of use for the conservation of the Lipschitz' constant. Still:
 $$\|Dh\|=\max(\|DR\|, \ \|DU\|)=|\beta| \leq  \gamma_0,$$
 hence $ \|Dh\| \|B^{-1}\|(1+\gamma_0)\leq \gamma_0 |\alpha|^{-1}(1+\gamma_0)<1$ since $|\alpha|\geq \cos \gamma_0$ and $\gamma_0$ is small enough. In particular, $P^{-1} \circ (\Psi_0(Y),Y)$ is a graph $(\Phi(Y),Y)$ which satisfies:
\begin{align*}
|\Phi(Y_1)-\Phi(Y_2)|&\leq \frac{\|A\| \gamma_0+ \|Dh\|(1+\gamma_0)}{\|B^{-1}\|^{-1}-\|Dh\|(1+\gamma_0)}|Y_1-Y_2|\\
                     &\leq \frac{|\alpha|\gamma_0+ \gamma_0(1+\gamma_0)}{|\alpha|-\gamma_0(1+\gamma_0)}|Y_1-Y_2|\leq C(\gamma_0)|Y_1-Y_2|.
\end{align*} 
 Finally, $(\Psi_0(Y),Y)$ is a graph for $|Y|\leq e^{-\chi_2n-9\gamma n}$ hence $(\Phi(Y),Y)$ is a graph for:
 $$ |Y|\leq \left(\|B^{-1}\|^{-1}- \|Dh\|(1+\gamma_0)\right)e^{-\chi_2n-9\gamma n} . $$
  We have that $(\|B^{-1}\|^{-1}- \|Dh\|(1+\gamma_0))e^{-\chi_2n-9\gamma n} \geq e^{-\chi_2n-10\gamma n}$ since
  $$\|B^{-1}\|^{-1}- \|Dh\|(1+\gamma_0)\geq  \alpha -\gamma_0(1+\gamma_0) \geq \frac{1}{2} $$
  for $\gamma_0$ small enough. \hfill $\Box$ \hfill \\
  
  We only consider the part $W' \subset W$ of points $(\Phi(Y),Y)$ with $|Y|\leq e^{-\chi_1 n -12 \gamma n}$. We consider the form:
  $$\beta':= \int_{W'} [D'_1(z) \cap B'] d\Hcal^2(z)$$
 where $D'_1(z)$ is the notation for the leaf of $\beta_1$ passing through $z\in W'$ (there is only one since $W$ is a graph over $D_1^{\bot}$), $\Hcal^2$ is the $2$-dimensional Hausdorff measure on $W'$ and $B':=B(0,e^{-\chi_1 n -12 \gamma n})$. Observe that $\beta'$ is still smooth since it can be seen as the pull-back of a smooth measure by a submersion.
We compare $\beta'$ and $\beta_1$:
 \begin{lemma}
 We have that $\beta'\geq \beta_1$ on $B'$.
 \end{lemma} 
 \noindent \emph{Proof.} Let $\Theta$ be a smooth positive $(1,1)$-form compactly supported in $B'$. We show that $\langle \beta', \Theta \rangle \geq \langle \beta_1, \Theta \rangle $. Consider the disk $D':= D_1^{\bot} \cap B'$ and let $Q:D'\to W$ be the map that sends $z$ to the unique point of intersection of $W$ with the leaf of 
 $\beta_1$ passing through $z$. In particular, $W'=Q(D')$.    
 As $W'$ is the set of $(\Phi(Y),Y)$ with $Y \in D'$, we have that $Q(Y)=(\Phi(Y),Y)$, in particular:
 $\|DQ(Y)\|=\|(\Phi'(Y),1)\|\geq 1 $ where $D$ is the complex differential ($Q$ is holomorphic). Recall that $D_1(z)$ is the leaf of $\beta_1$ passing through $z \in D_1^\bot$ so that $D'_1(Q(Y))=D_1(Y)$ with our notations. 
 Using that $Q$ is a bijection gives:
 \begin{align*}
 \int \beta'\wedge \Theta = \int_{W'} \left(\int_{D'_1(z)} \Theta\right) d\Hcal^2(z)\\
                          = \int_{Q(D')} \left(\int_{D_1(Q^{-1}z)} \Theta\right) d\Hcal^2(z).
 \end{align*}
 We can write $L(Y):= \int_{D_1(Y)} \Theta$ so that the previous equality is 
 $$\langle \beta',\Theta \rangle=\int_{Q(D')} L(Q^{-1}(z)) d\Hcal^2(z)  .$$ 
 The coarea formula (\cite{Fed}[p.258]) implies ($Q$ is a bijection):
 \begin{align*}
\langle \beta', \Theta \rangle  &= \int_{D'} L(Y) \|DQ(Y)\|  d\Hcal^2(Y)\\
                          &\geq \int_{D'} L(Y) d\Hcal^2(Y) = \int_{D'} \int_{D_1(Y)} \Theta  idY\wedge d\bar{Y}=\langle \beta_1, \Theta \rangle, 
 \end{align*} 
 where we use that since $D'$ is linear, the Hausdorff measure on $D'$ is just the standard Lebesgue measure. This is what we want. \hfill $\Box$ \hfill \\ 
 
 In particular, we have that (\ref{eqchi_1}) gives:
\begin{equation}\label{eqchi_2}
\Scal \wedge \beta'(B(0,e^{-\chi_1n-14\gamma n })) \geq \frac{C(\alpha_0,\gamma_0)}{2}e^{(-\chi_1n-15\gamma n)(\bar{d}_S+\gamma)}.
\end{equation} 
 We will push forward that inequality by $g^n_x$. We denote by $B$ the ball $B(0,e^{-\chi_1n-14\gamma n })$. Observe that $g^n_x$ is a diffeomorphism from $B'$ to $g^n_x(B')$. Indeed, for all $0\leq l \leq n-1$, $ g_x^l(B') \subset B(0,e^{-10 \gamma n}) $ by Lemma \ref{lemma3} and $g^{-1}_{\hat{f}^{n-l}(\hat{x})}$ is defined on $B(0, e^{-6\gamma n})$ (see the beginning of the proof  of Lemma \ref{Backward graph transform}). So, the fact that the map is meromorphic is indeed not an issue here.
 
  Let $0\leq \chi_B \leq 1$ be a smooth cut-off function equal to $1$ in the ball   $B$ and with support in $B'':=B(0, e^{-\chi_1 n -13.5 \gamma n})$ (so that $B \Subset B''\Subset B'$). By (\ref{eqchi_2}), the positive measure $\chi_B \Scal \wedge \beta'$ satisfies  $$\chi_B \Scal \wedge \beta' (B) \geq \frac{C(\alpha_0,\gamma_0)}{2}e^{(-\chi_1n-15\gamma n)(\bar{d}_S+\gamma)}.$$ 
 Pushing-forward implies that the measure $(g_x^n)_*(\chi_B \Scal \wedge \beta')$ has the same mass on $g_x^n(B)$:
 \begin{equation}\label{p33}
 (g^n_x)_*(\chi_B \Scal \wedge \beta') (g^n_x(B)) \geq \frac{C(\alpha_0,\gamma_0)}{2}e^{(-\chi_1n-15\gamma n)(\bar{d}_S+\gamma)}
 \end{equation}
 
  As $g^n_x$ is a diffeomorphism and $\beta'$ is smooth,
we have on $(g_x^n)(B'')$:
$$   (g_x^n)_*(\chi_B \Scal \wedge \beta')=(g_x^n)_*(\chi_B \Scal) \wedge  (g_x^n)_*(\beta') . $$ 
 \begin{proposition}\label{backtof}
In $g_x^n (B'')$, we have:
$$(g_x^n)_*(\beta')\leq e^{-2\chi_2 n + 18 \gamma n} \times (\tau_{f^n(x)}\circ C_\gamma(\hat{f}^n(\hat{x})))^*\omega.$$
\end{proposition}
 \noindent \emph{Proof.} Let $\Theta$ be a smooth positive $(1,1)$ form with compact support in $g_x^n (B'')$. Then the result will follow from:
$$ \langle (g_x^n)_*(\beta'), \Theta \rangle \leq  \langle e^{-2\chi_2 n + 18 \gamma n} \times (\tau_{f^n(x)}\circ C_\gamma(\hat{f}^n(\hat{x})))^*\omega, \Theta\rangle.$$
We have:
\begin{align*}
\int (g_x^n)_*(\beta') \wedge \Theta &=\int \beta' \wedge (g_x^n)^*(\Theta) \\
                                     &= \int_{W'} \int_{D'_1(z) \cap B'} (g_x^n)^*(\Theta) d\Hcal^2(z)\\
																		 &= \int_{W'} \int_{(g_x^n)(D'_1(z) \cap B')} \Theta d\Hcal^2(z)
\end{align*} 
where $W'$ is the part of $W$ of points $(\Phi(Y),Y)$ with $|Y|\leq e^{-\chi_1 n -12 \gamma n}$, and $W$ is a graph  above $D_1^\bot$ of Lipschitz' constant $\leq C(\gamma_0)$. Observe that:
\begin{align*}
\{z \in W', \ (g_x^n)(D_1'(z) \cap B') \cap  g_x^n (B'')  \neq \varnothing    \} &=\{z \in W',  \ D_1'(z) \cap B'' \neq \varnothing    \} \\
                                                                               &\subset  B(0, \sqrt{1+C(\gamma_0)^2}.e^{-\chi_1 n-13.5\gamma n}) \\
																																							 & \subset B'.
\end{align*} 
In particular:
\begin{align*}
\int (g_x^n)_*(\beta') \wedge \Theta = \int_{W' \cap B'} \int_{(g_x^n)(D'_1(z) \cap B')} \Theta d\Hcal^2(z).
\end{align*} 
We now apply the coarea formula: 
\begin{align}\label{coarea2}
  &\int_{W'\cap B'} \left(\int_{g^n_x(D'_1(z)\cap B')}\Theta\right) \|D({g^n_x}_{| W'\cap B'})(z)\|^2 d\Hcal^2(z) =\\
	&\int_{g^n_x(W'\cap B')} \left(\int_{g^n_x(D'_1((g^n_x)^{-1}(y))\cap B')}\Theta\right)  d\Hcal^2(y). \nonumber
\end{align}
Observe that the leaves $g^n_x(D'_1(z)\cap B')$ fill $g^n_x( B')$. We denote by $D_1^n(y)$ the leaf given by the image $g^n_x(D'_1(z)\cap B')$ that contains $y \in g^n_x( B')$ so that the previous term can be written:
$$
\int_{g^n_x(W'\cap B')} \left(\int_{D_1^n(y)}\Theta\right)  d\Hcal^2(y).
$$
Finally, $g^n_x(W'\cap B')\subset C_\gamma^{-1}(\hat{f}^n(\hat{x}))E_2(\hat{f}^n(\hat{x}))$ which is linear hence the Hausdorff measure $d\Hcal^2(y)$ can be replaced by $idy\wedge d\bar{y}$. We need the following lemma to complete the proof:
\begin{lemma}\label{lemma9}
For $z \in W'\cap B'$, we have:
$$e^{\chi_2 n + 2 \gamma n}\geq \|D({g^n_x}_{| W'\cap B'})(z)\| \geq e^{\chi_2 n - 2 \gamma n}  .$$
\end{lemma}
 \noindent \emph{Proof of the Lemma.}  Let $z \in W'\cap B'$ and $v \in T_zW'$, $\|v\|=1$ a vector tangent to $W'$ in $z$. We bound $\|Dg^n_x(z)v\|$. Recall that $g_x^n=g_{f^{n-1}(x)}\circ \dots  \circ g_x$, hence:
\begin{align} \label{etoile}
\|Dg^n_x(z)v\|&= \|Dg_{f^{n-1}(x)}(g^{n-1}_x(z)) \circ  Dg^{n-1}_x(z) v\|  \nonumber \\ 
              &= \frac{\|Dg_{f^{n-1}(x)}(g^{n-1}_x(z)) \circ  Dg^{n-1}_x(z) v\| }{\|Dg^{n-1}_x(z) v\|}\|Dg^{n-1}_x(z) v\|   \nonumber \\
              &= \frac{\|Dg_{f^{n-1}(x)}(g^{n-1}_x(z)) \circ  Dg^{n-1}_x(z) v\| }{\|Dg^{n-1}_x(z) v\|}\times   \nonumber \\
							& \ \frac{\|Dg_{f^{n-2}(x)}(g^{n-2}_x(z)) \circ  Dg^{n-2}_x(z) v\| }{\|Dg^{n-2}_x(z) v\|} \times \dots \times \frac{\|Dg_x(z)v \|}{\|v\|} \|v\|.
\end{align}
So it suffices to bound $\|Dg_{f^l(x)}(g^l_x(z)) u\|$ for $u$ a unit vector tangent to $g^l_x(W')$ at $g^l_x(z)$ for $l=0,\dots, n-1$. By construction, $g^l_x(W')$ is in the graph $(\Psi_l(Y), Y)$ with $Y\in C^{-1}_\gamma(\hat{f}^l(\hat{x}))E_2(\hat{f}^l(\hat{x}))$ and with $\mathrm{Lip}(\Psi_l)\leq \gamma_0$. In particular: 
$$ u = \frac{(\Psi'_l(Y),1)  }{\|(\Psi'_l(Y),1) \| }, \ \mathrm{where} \ (\Psi_l(Y),Y)=g_x^l(z).$$
We decompose:
$$\|Dg_{f^l(x)}(g^l_x(z)) u\|= \|Dg_{f^l(x)}(g^l_x(0)) u +Dg_{f^l(x)}(g^l_x(z)) u -    Dg_{f^l(x)}(g^l_x(0)) u        \| .$$
Hence:
\begin{align*}
&\|Dg_{f^l(x)}(g^l_x(0)) u \|+\|Dg_{f^l(x)}(g^l_x(z)) u -    Dg_{f^l(x)}(g^l_x(0)) u        \| \\
&\qquad \qquad \qquad \ \geq \|Dg_{f^l(x)}(g^l_x(z)) u\|\geq  \\
&\|Dg_{f^l(x)}(g^l_x(0)) u \| - \|Dg_{f^l(x)}(g^l_x(z)) u -    Dg_{f^l(x)}(g^l_x(0)) u        \| 
\end{align*}
We use the notations of Lemma \ref{Forward graph transform} and Corollary \ref{estimedh}:
$$g_{f^l(x)}(X,Y)= (AX + R(X,Y), B Y +U(X,Y))  $$
 with $A= A^1_\gamma(\hat{f}^l(\hat{x}))$, $ B= A^2_\gamma(\hat{f}^l(\hat{x}))$ and we write $Dh=(DR, DU)$. Hence, we have the computations:
\begin{align*}
&\|Dg_{f^l(x)}(0)u\|=\| (A,B)u\|= \frac{\|A \Psi'_l(Y)+B\|}{\|(\Psi'_l(Y),1)\|},\\
&\frac{\|B\| + \|A\| |\Psi'_l(Y)|}{\|(\Psi'_l(Y),1)\|}\geq\|Dg_{f^l(x)}(0)u\|\geq \frac{\|B\| - \|A\| |\Psi'_l(Y)|}{\|(\Psi'_l(Y),1)\|},\\
&e^{\chi_2+\gamma} + \gamma_0 e^{\chi_1+\gamma}  \geq\|Dg_{f^l(x)}(0)u\|\geq \frac{e^{\chi_2-\gamma} - \gamma_0 e^{\chi_1+\gamma}}{\sqrt{1+\gamma_0^2}} .
\end{align*}
Furthermore, since $g^l_x(z)\in g^l_x(B') \subset B(0,e^{-10 \gamma n})$ by Lemma \ref{lemma3} and Corollary \ref{estimedh}: $\|Dg_{f^l(x)}(g^l_x(z)) u -    Dg_{f^l(x)}(g^l_x(0)) u        \| \leq \| Dh(g^l_x)(z) \| \leq e^{-5\gamma n}$. Finally:
\begin{align*}
e^{\chi_2+\gamma} + \gamma_0 e^{\chi_1+\gamma} +e^{-5\gamma n} \geq\|Dg_{f^l(x)}(g^l_x(z)) u\|\geq \frac{e^{\chi_2-\gamma} - \gamma_0 e^{\chi_1+\gamma}}{\sqrt{1+\gamma_0^2}} -e^{-5\gamma n} .
\end{align*}
Recall that we chose $\gamma_0$ so it satisfies (\ref{gamma0}) so taking $n$ large enough gives:
\begin{align*}
e^{\chi_2+2\gamma} \geq \|Dg_{f^l(x)}(g^l_x(z)) u\| \geq e^{\chi_2-2\gamma}.
\end{align*}
The lemma follows. \hfill $\Box$ \hfill \\

\noindent \emph{Continuation of the proof of Proposition \ref{backtof}.} We apply the previous lemma to (\ref{coarea2}):
$$ \int_{W'\cap B'} \left(\int_{g^n_x(D'_1(z)\cap B')}\Theta\right)  d\Hcal^2(z) \leq e^{-2 \chi_2 n + 4\gamma n}	\int_{g^n_x(W'\cap B')} \left(\int_{D_1^n(y)}\Theta\right) idy \wedge d\bar{y}.$$
Consider the map $Q_1: g^n_x(B') \to C_\gamma^{-1}(\hat{f}^n(\hat{x})) E_2(\hat{f}^n(\hat{x}))$ defined as follows: for $y \in g^n_x(B') $ then $y \in D^1_n(y)$, by Lemma \ref{Forward graph transform}, it is a graph above  $ C_\gamma^{-1}(\hat{f}^n(\hat{x})) E_1(\hat{f}^n(\hat{x}))$ and we let $Q_1(y)$ be the unique projection of that graph with $C_\gamma^{-1}(\hat{f}^n(\hat{x})) E_2(\hat{f}^n(\hat{x}))$. The map $Q_1$ can also be defined as $g_x^n \circ Q_2 \circ (g_x^n)^{-1}$ where $Q_2$ is the projection from $B'$ to $W'$ that sends $z' \in D'_1(z)$ to $z$.

 The current $\int_{g^n_x(W'\cap B')} [D_1^n(y)] idy \wedge d\bar{y}$ can be written as the pull back of $(idy\wedge d\bar{y})_{|g^n_x(W'\cap B')}$ by the submersion $Q_1$. Hence:
$$ \int (g_x^n)_*(\beta') \wedge \Theta \leq e^{-2 \chi_2 n + 4\gamma n}	\int_{g^n_x(W'\cap B')} \Theta \wedge  Q_1^*((idy\wedge d\bar{y})_{|g^n_x(W'\cap B')}).$$ 
We need the lemma:
\begin{lemma}\label{lemma10}
We have that in $g^n_x(B')$:
$$Q_1^*((idy\wedge d\bar{y})_{|g^n_x(W')}) \leq e^{11 \gamma n} \beta_0,$$
where $\beta_0$ denotes the standard Kähler form in $\Cc^2$. 
\end{lemma}
\noindent \emph{Proof of the Lemma.} By the chain rule:
$$DQ_1(z)= Dg^n_x(Q_2((g^n_x)^{-1}(z))) \circ DQ_2((g^n_x)^{-1}(z)) \circ D(g^n_x)^{-1}(z).$$  
Let $v$ be a unit vector with $DQ_1(z)v\neq 0$. We have:
\begin{align*}
\|DQ_1(z)v\|= \frac{\|Dg^n_x(Q_2((g^n_x)^{-1}(z)))(DQ_2((g^n_x)^{-1}(z))( D(g^n_x)^{-1}(z))) v    \|}{ \| DQ_2((g^n_x)^{-1}(z))( D(g^n_x)^{-1}(z))v\|  }\\
\times   \frac{\| DQ_2((g^n_x)^{-1}(z))( D(g^n_x)^{-1}(z))v\|}{\|D(g^n_x)^{-1}(z)v\|}
\times \|D(g^n_x)^{-1}(z)v\|.
\end{align*}
Since $Q_2((g^n_x)^{-1}(z)) \in W'$ for $z \in g^n_x(B')$, using Lemma \ref{lemma9} gives:
\begin{align*}
 \frac{\|Dg^n_x(Q_2((g^n_x)^{-1}(z)))(DQ_2((g^n_x)^{-1}(z))( D(g^n_x)^{-1}(z))) v    \|}{ \| DQ_2((g^n_x)^{-1}(z))( D(g^n_x)^{-1}(z))v\|  }\leq e^{\chi_2 n +2\gamma n}.
\end{align*}
Since $W$ is a graph above $D^\bot_1$ which is $C(\gamma_0)$-Lipschitz (Lemma \ref{lemma6}), we have:
\begin{align*}
\frac{\| DQ_2((g^n_x)^{-1}(z))( D(g^n_x)^{-1}(z))v\|}{\|D(g^n_x)^{-1}(z)v\|}\leq C(\gamma_0).
\end{align*}
Finally, we have $(g^l_x)^{-1}(z) \in B(0,e^{-10 \gamma n}) $ for $l=0,\dots, n$ by Lemma \ref{lemma3}   and  $(g^n_x)^{-1}= (g^{-1}_{\hat{f}(\hat{x})}) \circ \dots \circ (g^{-1}_{\hat{f}^n(\hat{x})}) $. We have for $l=0,\dots,n$, using the notations of formula (\ref{estimeDh2}) and Corollary \ref{estimedh2}: 
\begin{align*} 
 \|Dg^{-1}_{\hat{f}^l(\hat{x})}(w) \| &\leq  \|Dg^{-1}_{\hat{f}^l(\hat{x})}(w)-Dg^{-1}_{\hat{f}^l(\hat{x})}(0) \|+ \|Dg^{-1}_{\hat{f}^l(\hat{x})}(0) \| \\
                                       &\leq \|Dh(w)\| + \|Dg^{-1}_{\hat{f}^l(\hat{x})}(0) \|\\
                                       &\leq e^{-5 \gamma n} +e^{-\chi_2  + \gamma}\leq e^{-\chi_2  +2 \gamma}
\end{align*}
for $w \in B(0,e^{-10 \gamma n})$. Hence $\|D(g^n_x)^{-1}(z) v \|\leq \exp(-\chi_2 n + 2\gamma n) $ for $z \in g^n_x(B')$. Finally $ \|DQ_1(z)v\|\leq e^{5 \gamma n}$ for $z \in g^n_x(B')$. So we have $Q_1^*(idy\wedge d\bar{y})=idQ_1\wedge d\bar{Q_1} \leq e^{11 \gamma n}\beta_0$ in $g^n_x(B')$. \hfill $\Box$ \hfill \\

\noindent \emph{End of the proof of Proposition \ref{backtof}.} Using the above lemma, we deduce:
$$ \int (g^n_x)_* \beta' \wedge \Theta \leq e^{-2\chi_2  n +4 \gamma n} e^{11\gamma n} \int \Theta \wedge \beta_0. $$ 
This means that $(g^n_x)_*\beta' \leq e^{-2\chi_2 n +15 \gamma n}\beta_0$ in $g^n_x(B'')$ ($\Theta$ is arbitrary with support in $g^n_x(B'')$). Since 
$$\beta_0\leq C \frac{1}{\alpha_0^2} e^{2 \gamma n} (\tau_{f^n(x)} \circ C_\gamma(\hat{f}^n(\hat{x})))^* \omega$$
the proposition follows. \hfill $\Box$ \hfill \\

We finish the proof of  Theorem \ref{Theorem1} in the case where $\chi_1 > \chi_2$. We have by (\ref{p33}) and Proposition \ref{backtof}:
 \begin{align*}
 e^{-2\chi_2 n +18 \gamma n}  \left((g^n_x)_*(\chi_B \Scal) \wedge (\tau_{f^n(x)} \circ C_\gamma(\hat{f}^n(\hat{x})))^*\omega \right)   (g^n_x(B))  \\
 \geq \frac{C(\alpha_0,\gamma_0)}{2}e^{(-\chi_1n-15\gamma n)(\bar{d}_S+\gamma)}.
 \end{align*}So:
 \begin{align*}
\left((g^n_x)_*(\chi_B \Scal) \wedge (\tau_{f^n(x)} \circ C_\gamma(\hat{f}^n(\hat{x})))^* \omega \right)    (g^n_x(B)) \geq \\ \frac{C(\alpha_0,\gamma_0)}{2}e^{(-\chi_1n-15\gamma n)(\bar{d}_S+\gamma)+2\chi_2 n -18\gamma n}.
 \end{align*}
Using that $ g^n_x=C_\gamma^{-1}(\hat{f}^n(\hat{x})) \circ \tau_{f^n(x)}^{-1} \circ f^n \circ \tau_{x} \circ C_\gamma(\hat{x})$ and  $\Scal = (\tau_x \circ C_\gamma(\hat{x}))^*(S)$ gives:
 \begin{align*}
& \left( (g^n_x)_*(\chi_B \Scal) \wedge (\tau_{f^n(x)} \circ C_\gamma(\hat{f}^n(\hat{x})))^* \omega\right)   (g^n_x(B)) =\\
& \left( (f^n\circ \tau_x\circ C_\gamma(\hat{x}))_*(\chi_B \Scal) \wedge \omega\right)  ( f^n\circ \tau_x\circ C_\gamma(\hat{x})(B)) =\\
&\left( (f^n)_*(\chi_B\circ (\tau_x\circ C_\gamma(\hat{x}))^{-1} S) \wedge \omega\right)  ( f^n\circ \tau_x\circ C_\gamma(\hat{x})(B)).
   \end{align*}
Hence:
 \begin{align}\label{utilepourth2}
(f^n)_*(\chi_B\circ (\tau_x\circ C_\gamma(\hat{x}))^{-1} S) \wedge \omega( (f^n\circ \tau_x\circ C_\gamma(\hat{x}))(B)) \geq  \nonumber \\ \frac{C(\alpha_0,\gamma_0)}{2}e^{(-\chi_1n-15\gamma n)(\bar{d}_S+\gamma)+2\chi_2 n -18\gamma n}.
 \end{align}
Now, let $1_A$ denote the indicator function of a set $A$, then
$$\chi_B\circ (\tau_x\circ C_\gamma(\hat{x}))^{-1}\leq 1_{ (\tau_x\circ C_\gamma(\hat{x}))(B'')} \leq 1_{B(x,e^{-\chi_1 n-13 \gamma n})} \leq 1_{B_n(x,\frac{\delta}{2})},  $$ 
by Lemma \ref{lemma3}. Now using the results of \cite{DS4}, we know that there exists a positive closed current $S'$ of mass $\leq C \|S\|$ where $C$ is a constant that does not depend on $S$ such that $S'$ is a limit of smooth positive closed current (in the projective case, one can take $S=S'$). Let $S_m$ be such a sequence. Then we have:
$$\lim_{m\to \infty} \int (f^n)_*( 1_{B_n (x,\frac{\delta}{2})} S_m )\wedge \omega  \geq  \frac{C(\alpha_0,\gamma_0)}{2}e^{(-\chi_1n-15\gamma n)(\bar{d}_S+\gamma) +2\chi_2 n-18\gamma n}.$$
We sum on all $x=x_i$ where the $x_i$ are the points of the $(n, \delta)$-separated set so that the $B_n(x_i,\delta/2)$ are disjoint:
$$\lim_{m\to \infty} \sum_i \int (f^n)_*( 1_{B_n (x_i,\frac{\delta}{2})} S_m )\wedge \omega \geq  N \frac{C(\alpha_0,\gamma_0)}{2}e^{(-\chi_1n-15\gamma n)(\bar{d}_S+\gamma) +2\chi_2 n -18\gamma n}.$$
One knows that there exists a constant $\alpha_X$ such that for any smooth positive closed current $S_m$  of mass $1$, $\alpha_X \omega - S_m$ is positive up to a $dd^c$-closed current (see e.g. \cite{DS9}[remarque 3]), in particular:
 $$C \alpha_X \lambda_1(f^n) \geq  \int (f^n)_*( S_m )\wedge \omega$$ 
by definition of $\lambda_1(f^n)$. Replacing the constant $C \alpha_X$ by $C$ for simplicity gives:
$$C   \lambda_1(f^n) \geq  \int (f^n)_*( S_m )\wedge \omega  \geq  \sum_i \int (f^n)_*( 1_{B_n (x_i,\frac{\delta}{2})} S_m )\wedge \omega. $$
  Recall that $N \geq \frac{1}{4} e^{h_\nu(f)n-\gamma n}$ so we have:
$$  C  \lambda_1(f^n) \geq \frac{C(\alpha_0,\gamma_0)}{8} e^{h_\nu(f)n-\gamma n} e^{(-\chi_1n-15\gamma n)(\bar{d}_S+\gamma) +2\chi_2 n -18\gamma n }.  $$ 
We take the logarithm, we divide by $n$:
$$\frac{\log \lambda_1(f^n)}{n} \geq \frac{1}{n}\log ( \frac{C(\alpha_0,\gamma_0)}{8C}) +h_\nu(f)-19 \gamma +2 \chi_2  -(\chi_1 +15 \gamma)(\bar{d}_S+\gamma),$$
we let $n\to \infty$ then $\gamma \to 0$:
$$\log d_1 \geq h_\nu(f)+2 \chi_2-\chi_1 \bar{d}_S,$$
which concludes the proof in the case where $\chi_1 > \chi_2$.

\subsection{Case where $\chi_1 = \chi_2=\chi>0$.} 
We follow the same approach than in the previous case. This time we do not have any privileged direction but the dilatation is almost uniform.  We still denote by $x$ one of the $x_i$ in the $(n, \delta)$-separated set. In particular, we have:
$$ S\wedge \omega (B(x,e^{-\chi n -15 \gamma n})) \geq e^{(-\chi n - 15 \gamma n)(\bar{d}_s+\gamma)}. $$
We write that inequality in the local chart of $\Cc^2$. The standard Kähler form $\beta$ of $\Cc^2$ is written 
$\beta_1 +\beta_2$ where $\beta_j = i dz_j\wedge d\bar{z}_j$ ($j=1,2$). 

We have in a ball $B(x, C \varepsilon_0) \supset B(x,e^{- \gamma n})$:
$$ (\tau_x \circ C_\gamma(\hat{x}))_* \beta \geq C(\alpha_0) \omega $$
for some constant $C(\alpha_0)$. Again, we denote in what follows $\Scal := (\tau_x \circ C_\gamma(\hat{x}))^*(S)$. 
We consider the balls $B:=B(0,e^{-\chi n-14\gamma n })$, $B':=B(0,e^{-\chi n -12 \gamma n})$ and $B'':=B(0, e^{-\chi  n -13.5n})$.
Then, we proceed exactly as in the proof of formula (\ref{eqchi_1}) and we obtain:
\begin{equation}
\label{eqchi_12}
\Scal \wedge \beta_1(B) \geq  \frac{C(\alpha_0)}{2} e^{(-\chi n-15 \gamma n)(\bar{d}_S+\gamma)}.
\end{equation}
We push that inequality forward by $g_x^n$, it is a diffeomorphism from $B'$
to $g_x^n(B')$ (see the paragraph after (\ref{eqchi_2})). We consider again  $0\leq \chi_B \leq 1$, a smooth cut-off function equal to $1$ in the ball $B$ and with support in $B''$. Then $\chi_B \Scal \wedge \beta_1$ is a positive measure of mass $\geq  \frac{C(\alpha_0)}{2} e^{(-\chi  n-15 \gamma n)(\bar{d}_S+\gamma)} $ and $(g_x^n)_*(\chi_B \Scal \wedge \beta_1)$ has the same mass.  Let $D:=(z_2=0)$, then we can decompose $\beta_1$ as :
$$\beta_1:=\int_{D }[\{(a,b)\in \Cc^2, \ a=z_1\}] idz_1\wedge d\bar{z}_1 .$$
Let $\beta'$ denote:
$$\beta':=\int_{D }[\{(a,b)\in \Cc^2, \ a=z_1\}\cap B'] idz_1\wedge d\bar{z}_1.$$
 Now 
\begin{align*}
(g_x^n)_*(\chi_B \Scal \wedge \beta_1)&=(g_x^n)_*(\chi_B \Scal \wedge \beta') \\
                                      &=(g_x^n)_*(\chi_B \Scal) \wedge  (g_x^n)_*(\beta'). 
\end{align*}
 \begin{proposition}\label{backtof2}
In $g_x^n (B'')$, we have:
$$(g_x^n)_*(\beta')\leq e^{-2\chi n  + 18 \gamma n} \times(\tau_{f^n(x)}\circ C_\gamma(\hat{f}^n(\hat{x})))^*\omega.$$
\end{proposition}
 \noindent \emph{Proof.} Let $\Theta$ be a smooth positive $(1,1)$ form with compact support in $g_x^n (B'')$. Then the result will follow from:
$$ \langle (g_x^n)_*(\beta'), \Theta \rangle \leq  \langle e^{-2\chi n  + 18 \gamma n} .(\tau_{f^n(x)}\circ C_\gamma(\hat{f}^n(\hat{x})))^*\omega, \Theta\rangle.$$
We have:
\begin{align*}
\int (g_x^n)_*(\beta') \wedge \Theta &=\int \beta' \wedge (g_x^n)^*(\Theta) \\
                                     &= \int_{D} \int_{\{(a,b), \ a=z_1\}\cap B'} (g_x^n)^*(\Theta) idz_1\wedge d\bar{z}_1\\
																		 &= \int_{D} \int_{g_x^n(\{(a,b), \ a=z_1\}\cap B')} \Theta idz_1\wedge d\bar{z}_1\\
						                         &= \int_{D\cap B'} \int_{g_x^n(\{(a,b), \ a=z_1\}\cap B')} \Theta d\Hcal^2(z_1),						                        							 
\end{align*}
where we use for the last equality that:
$$\{z_1, \ g_x^n(\{(a,b), \ a=z_1\})\cap g_x^n(B') \neq \varnothing \} = \{z_1, (\{(a,b), \ a=z_1\})\cap B' \neq \varnothing \}, $$
and where $\Hcal^2(z_1)$ is the Hausdorff measure on $D$. 
We now apply the coarea formula: 
\begin{align}\label{coarea22}
  &\quad \int_{D\cap B'} \left(\int_{g_x^n(\{(a,b), \ a=z_1\}\cap B')} \Theta\right)  \| (Dg^n_x)_{|D\cap B'}\|^2 d\Hcal^2(z_1)\\
	&= \int_{g^n_x(D\cap B')} \left(\int_{g_x^n(\{(a,b), \ a=(g^n_x)^{-1}(y) \} \cap B' )} \Theta \right)   d\Hcal^2(y). \nonumber
\end{align}
As before, the leaves $g_x^n(\{(a,b), \ a=(g^n_x)^{-1}(y) \} \cap B')$ fill $g^n_x( B')$ and we denote by $D_1^n(y)$ the leaf given by the image $g_x^n(\{(a,b), \ a=(g^n_x)^{-1}(y) \} \cap B')$ that contains $y \in g^n_x( B')$. So the previous term can be written:
$$\int_{g^n_x(D\cap B')} \left(\int_{D_1^n(y)}\Theta\right)  d\Hcal^2(y).$$
\begin{lemma}\label{lemma11}
For $z \in D\cap B'$, we have:
$$e^{\chi  n + 2 \gamma n}\geq \|D({g^n_x}_{| D\cap B'})(z)\| \geq e^{\chi  n - 2 \gamma n}  .$$
\end{lemma}
 \noindent \emph{Proof of the Lemma.}  Let $z \in D\cap B'$ and $v \in T_z(D\cap B')$, $\|v\|=1$ a vector tangent to $W'$ in $z$. We proceed as in the proof of Lemma \ref{lemma9}: we write the formula (\ref{etoile}) so it suffices to bound $\|Dg^{l}_{f^l(x)}(g^l_x(z)) u\|$ for $u$ a unit vector tangent to $g^l_x(D)$ at $g^l_x(z)$ for $l=0,\dots, n-1$. We write:
\begin{align*}
&\|Dg_{f^l(x)}(g^l_x(0)) u \|+\|Dg_{f^l(x)}(g^l_x(z)) u -    Dg_{f^l(x)}(g^l_x(0)) u        \| \\
&\qquad \qquad \qquad \ \geq \|Dg_{f^l(x)}(g^l_x(z)) u\|\geq  \\
&\|Dg_{f^l(x)}(g^l_x(0)) u \| - \|Dg_{f^l(x)}(g^l_x(z)) u -    Dg_{f^l(x)}(g^l_x(0)) u        \|. 
\end{align*}
Osedelec's Theorem implies:
\begin{align*}
 e^{\chi-\gamma}  \leq\|Dg_{f^l(x)}(0)u\|\leq e^{\chi+\gamma} 
\end{align*}
Furthermore, using the notations and results of Lemma \ref{lemma3} and Corollary \ref{estimedh}:
\begin{align*}
\|Dg_{f^l(x)}(g^l_x(z)) u -    Dg_{f^l(x)}(g^l_x(0)) u        \| \leq \| Dh(g^l_x)(z) \| \leq e^{-5\gamma n},
\end{align*}
since $g^l_x(z)\in g^l_x(B') \subset B(0,e^{-10 \gamma n})$. Finally:
\begin{align*}
e^{\chi+2\gamma} \geq \|Dg_{f^l(x)}(g^l_x(z)) u\| \geq e^{\chi-2\gamma}.
\end{align*}
The lemma follows. \hfill $\Box$ \hfill \\

Applying the previous lemma gives:
$$
\int_{D\cap B'} \left(\int_{g^n_x(\{(a,b), \ a=z_1\} \cap B')}\Theta\right) idz_1\wedge d\bar{z}_1 \leq e^{-2\chi n + 4\gamma n}\int_{g^n_x(D\cap B')} \left(\int_{D_1^n(y)}\Theta\right)  d\Hcal^2(y).$$
Consider the map $Q_1: g^n_x(B') \to g^n_x(D\cap B')$ defined by: for $y \in g^n_x(B') $ then $y \in g^n_x(\{(a,b), \ a=z_1\} \cap B')$ for a unique $z_1$, 
that set intersects $g^n_x(D\cap B')$ at a unique point that we denote by $Q_1(y)$. The map $Q_1$ can also be defined as $g_x^n \circ Q_2 \circ (g_x^n)^{-1}$ where $Q_2$ is the orthogonal projection from $B'$ to $D$.
The current $\int_{g^n_x(D\cap B')} [D_1^n(y)] idy \wedge d\bar{y}$ can be written as $Q_1^*((idy\wedge d\bar{y})_{|g^n_x(D \cap B')})$. Hence:
$$ \int (g_x^n)_*(\beta') \wedge \Theta \leq e^{-2 \chi n + 4\gamma n}	\int_{g^n_x(B')} \Theta \wedge  Q_1^*((idy\wedge d\bar{y})_{|g^n_x(D\cap B')}).$$ 
\begin{lemma}\label{lemma102}
We have that in $g^n_x(B')$:
$$Q_1^*((idy\wedge d\bar{y})_{|g^n_x(D\cap B')}) \leq e^{11 \gamma n} \beta_0,$$
where $\beta_0$ denotes the standard Kähler form in $\Cc^2$. 
\end{lemma}
The proof is exactly the same than for Lemma \ref{lemma10} so we skip it.  We finish the proof of Proposition \ref{backtof2}, we have:
\begin{align*}
\int (g_x^n)_*(\beta') \wedge \Theta &\leq e^{-2 \chi n + 4\gamma n}	\int_{g^n_x( B')} \Theta \wedge  Q_1^*((idy\wedge d\bar{y})_{|g^n_x(D\cap B')})\\
                                     & \leq e^{-2 \chi n +15 \gamma n } \int_{g^n_x( B')} \Theta \wedge \beta_0.
\end{align*}
That means $(g_x^n)_*(\beta') \leq e^{-2 \chi n +15 \gamma n } \beta_0 $ in $g^n_x(B'')$. Since
$$\beta_0\leq C \frac{1}{\alpha_0^2} e^{2 \gamma n} (\tau_{f^n(x)} \circ C_\gamma(\hat{f}^n(\hat{x})))^* \omega$$
the proposition follows. \hfill $\Box$ \hfill \\

We finish the proof in the case $\chi_1=\chi_2=\chi$. We have that from (\ref{eqchi_12}):
\begin{align*}
(g_x^n)_*(\chi_B \Scal \wedge \beta_1)(g_x^n(B))  &\geq \frac{C(\alpha_0)}{2} e^{(- \chi n - 15 \gamma n )(\bar{d}_S+\gamma)} \\
(g_x^n)_*(\chi_B \Scal) \wedge (g_x^n)_*(\beta')(g_x^n(B)) &\geq \frac{C(\alpha_0)}{2} e^{(- \chi n - 15 \gamma n )(\bar{d}_S+\gamma)}
\end{align*}
So the above proposition gives:
\begin{align*}
 (g_x^n)_*(\chi_B \Scal) \wedge (\tau_{f^n(x)} \circ C_\gamma(\hat{f}^n(\hat{x})))^* \omega (g_x^n(B)) \geq \\
 \qquad  \qquad \frac{C(\alpha_0)}{2} e^{2 \chi n -18 \gamma n - ( \chi n + 15 \gamma n )(\bar{d}_S+\gamma)}.
\end{align*}
We deduce:
 \begin{align}\label{utilepourth22}
(f^n)_*(\chi_B\circ (\tau_x\circ C_\gamma(\hat{x}))^{-1} S) \wedge \omega( (f^n\circ \tau_x\circ C_\gamma(\hat{x}))(B)) \geq  \nonumber \\
\frac{C(\alpha_0)}{2} e^{(- \chi n - 15 \gamma n )(\bar{d}_S+\gamma) + 2\chi n - 18 \gamma n}.
   \end{align}
Again:
$$\chi_B\circ (\tau_x\circ C_\gamma(\hat{x}))^{-1}\leq 1_{ (\tau_x\circ C_\gamma(\hat{x}))(B'')} \leq 1_{B(x,e^{-\chi  n-13 \gamma n})} \leq 1_{B_n(x,\frac{\delta}{2})},  $$ 
by Lemma \ref{lemma3}. The rest of the proof is the same than in the case where $\chi_1 > \chi_2$.

\section{Proof of Theorem \ref{Theorem2}}
We modify the proof of Theorem \ref{Theorem1}. Let $\gamma>0$ small compared to the Lyapunov exponents. We make stronger uniformizations. First, let:
\begin{align*}
 \Lambda_{k_0}:= \Big\{x\in \Lambda, \ \forall n \geq k_0, \ \frac{\log(S\wedge \omega(B(x,e^{-\chi_1n- 15 \gamma n})))}{\log e^{-\chi_1n- 15 \gamma n}} \leq \bar{d}_S +\gamma   \\ 
\text{and} \      \frac{\log(\nu(B(x,2e^{-\chi_1n- 12 \gamma n})))}{\log 2e^{-\chi_1n- 12 \gamma n}} \geq \underline{d}_\nu -\gamma     \Big\}.
\end{align*}
By (\ref{hypotheseS}) and the hypotheses of Theorem \ref{Theorem2}, $\bigcup_{k_0 \geq 0}\Lambda_{k_0}= \Lambda$ and $\nu(\Lambda)=1$. So we choose $k_0$ large enough such that $\nu(\Lambda_{k_0}) \geq 3/4$. 

Consider now :
$$X_{\delta, n}:= \left\{ x, \ \nu(B_n(x,5\delta)) \leq e^{-h_{\nu}(f)n +\gamma n}  \ \mathrm{and} \  \nu(B_n(x,\frac{\delta}{2})) \geq e^{-h_{\nu}(f)n -\gamma n}  \right\}. $$
Brin-Katok's formula gives that for $\delta$ small enough:
\begin{align*}
\frac{4}{5} &\leq \nu \big(\big\{x, \ \underline{\lim} -\frac{1}{n} \log \nu(B_n(x,5\delta))\geq h_\nu(f) -\frac{\gamma}{2}  \\
& \  \mathrm{and} \ \overline{\lim} -\frac{1}{n} \log \nu(B_n(x,\delta/2))\leq h_\nu(f) + \frac{\gamma}{2}       \big\}\big) \\
            &\leq \nu \left(  \bigcup_{n_0} \bigcap_{n \geq n_0} X_{\delta,n}   \right).
\end{align*} 
In particular, we choose $n_0$ large enough so that:
$$\nu\left(\bigcap_{n\geq n_0} X_{\delta,n_0}\right) \geq \frac{3}{4}.$$
We consider again:
$$\widehat{Y}_{\alpha_0}:=\left\{ \hat{x}\in \widehat{Y}, \ \alpha_0 \leq \| C^{\pm1}_\gamma (\hat{x})\|\leq \frac{1}{\alpha_0} \ \mathrm{and} \ V(\hat{x}) \geq \alpha_0 \right\}$$
where we have $d(x_n,\Acal)\geq V(\hat{x})e^{-\frac{\gamma|n|}{p}}$ for all $n\in \Zz$. Take $\alpha_0$ small enough so that $\hat{\nu}(\widehat{Y}_{\alpha_0})\geq 3/4$. 

Let $A':= \Pi(\widehat{Y}_{\alpha_0}) \cap \left(\bigcap_{n\geq n_0} X_{\delta,n_0}\right) \cap \Lambda_{k_0}$. Since $\Pi_*(\hat{\nu})=\nu$, we have that $\nu(\Pi(\widehat{Y}_{\alpha_0})) \geq \hat{\nu}(\widehat{Y}_{\alpha_0}) \geq 3/4$, hence $\nu(A')\geq 1/4$.

In what follows, we take $n \geq \max(n_0, k_0)$. If $x\in A'$, we have that $\nu(B_n(x,\delta))\leq e^{-h_{\nu}(f) n+\gamma n}$. Hence we can find a $(n,\delta)$ separated set $(x_1,\dots, x_{N_1})$ in $A'$ with $N_1 \geq \frac{1}{4} e^{h_{\nu}(f)n-\gamma n}$.
 Take that set maximal so that $A'\subset \cup_{i\leq N_1} B_n(x_i, \delta)$. Let $N_2$ be the number of $i$ such that $\nu(B_n(x_i,\delta)\cap A') \geq e^{-h_\nu(f) n 
-2 \gamma n}$. We bound $N_2$ from below:
\begin{align*} 
\frac{1}{4}  &\leq \nu(A') \leq \nu (\cup_{i\leq N_1} B_n(x_i, \delta) \cap A') \leq \sum _{i\leq N_1} \nu(B_n(x_i, \delta) \cap A') \\
            &\leq (N_1-N_2)e^{- h_\nu(f) n -2 \gamma n}+ N_2 e^{- h_\nu(f) n + \gamma n}.
 \end{align*}
Since the balls $B_n(x_i, \delta/2)$ are disjoint we have:
$$1 \geq \nu(\cup_{i\leq N_1} B_n(x_i, \delta/2))= \sum_{i\leq N_1} \nu( B_n(x_i, \delta/2))\geq N_1 e^{- h_\nu(f) n - \gamma n}.$$
Hence $ N_1-N_2 \leq N_1 \leq e^{ h_\nu(f) n + \gamma n}$ so: 
$$ \frac{1}{4}  \leq e^{ - \gamma n}+ N_2 e^{- h_\nu(f) n + \gamma n}.$$
This gives the bound
 $$N_2 \geq \frac{e^{ h_\nu(f) n - \gamma n}}{8}. $$
We simply denote the points satisfying $\nu(B_n(x_i,\delta)\cap A') \geq e^{-h_\nu(f) n 
-2 \gamma n}$ by $x_1$, \dots, $x_{N_2}$.
\begin{lemma}
We can extract from those $N_2$ points a $(n,4\delta)$ separated set of cardinality $N_3 \geq \frac{e^{ h_\nu(f) n - 3\gamma n}}{8}$.
\end{lemma}
\noindent \emph{Proof.} For each $i\leq N_2$, let $m_i$ be the maximal number of $x_j$ wih $x_j \in B_n(x_i, 4 \delta)$. Let $m:=\max m_i$  and $x$ be $x_i$ such that $m_i=m$. Let $x'_1, \dots, x'_m$ be the points in $B_n(x, 4 \delta)$. Then $\cup_{j=1}^m B_n(x'_j,\delta/2)\subset B_n(x, 5 \delta)$. As the balls $B_n(x'_j, \delta/2)$ are disjoint:
$$m e^{- h_\nu(f) n - \gamma n} \leq \sum_{j\leq m}  \nu(B_n(x'_j,\delta/2))= \nu(\cup_{j\leq m}  B_n(x'_j,\delta/2)) \leq \nu( B_n(x,5\delta)) \leq e^{- h_\nu(f) n + \gamma n} . $$
Hence $m \leq e^{2 \gamma n}$ and the lemma follows.  \hfill $\Box$ \hfill \\

We simply denote those points $x_1,\dots, x_{N_3}$. We fix $x$ one of those $x_i$. By construction, $\nu (B_n(x, \delta) \cap A') \geq e^{- h_\nu(f) n  -2 \gamma n}$ . For $y \in B_n(x, \delta) \cap A'$, we have:
$$ \nu( B(y, 2e^{-\chi_1 n -12 \gamma n}))\leq 2^{\underline{d}_\nu -\gamma } e^{(- \chi_1 n -12 \gamma n )(\underline{d}_\nu -\gamma ) }   .$$ 
We can thus find a $2e^{-\chi_1 n -12 \gamma n}$ separated set in $B_n(x, \delta)\cap A'$ whose cardinality $L$ satisfies:
$$ L \geq \frac{ e^{- h_\nu(f) n  -2 \gamma n}}{2^{\underline{d}_\nu -\gamma } e^{(- \chi_1 n -12 \gamma n )(\underline{d}_\nu -\gamma ) }}.$$
Let $y_1, \dots, y_L$ be those points (the balls $B(y, e^{-\chi_1 n -12 \gamma n})$ are disjoint). Let $y$ be one of those $y_j$. We now follow the proof of Theorem 1. Observe that $y \in A' \subset A$ hence it satisfies the same estimates than the point $x$ of the previous section. In particular, recall that we denoted $B=B(0,e^{-\chi_1 n - 14 \gamma n})$ and  $0\leq \chi_B \leq 1$ be a smooth cut-off function equal to $1$ in the ball $B$ and with support in $B'':=B(0, e^{-\chi_1 n -13.5n})$. We obtained (see (\ref{utilepourth2}) and (\ref{utilepourth22})):
 \begin{align*}
(f^n)_*(\chi_B\circ (\tau_y\circ C_\gamma(\hat{y}))^{-1} S) \wedge \omega( (f^n\circ \tau_y\circ C_\gamma(\hat{y}))(B)) \geq  \nonumber \\ \frac{C(\alpha_0,\gamma_0)}{2}e^{(-\chi_1n-15\gamma n)(\bar{d}_S+\gamma)+2\chi_2 n -18\gamma n}.
 \end{align*}
and $\chi_B\circ (\tau_y\circ C_\gamma(\hat{y}))^{-1} \leq 1_{B(y,e^{-\chi_1 n- 13 \gamma n})}$. So that:
$$ \int (f^n)_*( 1_{B(y, e^{-\chi_1 n- 13 \gamma n})} S ) \wedge \omega \geq  \frac{C(\alpha_0,\gamma_0)}{2}e^{(-\chi_1n-15\gamma n)(\bar{d}_S+\gamma) +2\chi_2 n -18\gamma n} .   $$
For a given $x$ (among the $N_3$) we have $L$ such $y$. As the balls $B(y_j,e^{-\chi_1 n- 13 \gamma n})$ are disjoint we have:
$$ \int (f^n)_*( \sum_{j=1}^{L} 1_{B(y_j, e^{-\chi_1 n- 13 \gamma n})} S ) \wedge \omega \geq  L\frac{C(\alpha_0,\gamma_0)}{2}e^{(-\chi_1n-15\gamma n)(\bar{d}_S+\gamma) +2\chi_2 n -18\gamma n} .   $$
By Lemma \ref{lemma3}, for each $j$ we have $B(y_j,e^{-\chi_1 n- 13 \gamma n}) \subset B_n(y_j, \delta /2 )$. As $y_j \in B_n(x, \delta)$, we deduce:
$$\sum_{j=1}^L 1_{B(y_j, e^{-\chi_1 n- 13 \gamma n})} \leq 1_{B_n(x , 2 \delta)}.$$
So, using the same approximation's argument of $S$ than in the proof of Theorem \ref{Theorem1}:
\begin{align*}
 \lim_{m\to  \infty}  \int (f^n)_*( 1_{B_n(x , 2 \delta)} S_m ) \wedge \omega \geq   \frac{ e^{- h_\nu(f) n  -2 \gamma n}}{2^{\underline{d}_\nu -\gamma }  e^{(- \chi_1 n -12 \gamma n )(\underline{d}_\nu -\gamma ) }} \times   \\
\frac{C(\alpha_0,\gamma_0)}{2}e^{(-\chi_1n-15\gamma n)(\bar{d}_S+\gamma) +2\chi_2 n -18\gamma n} .
\end{align*}
Finally, the points $x_1,\dots , x_{N_3}$ are $(n,4 \delta)$ separated so the balls $B(x_i, 2 \delta)$ are disjoint hence:
\begin{align*}
 \lim_m \int (f^n)_*( \sum_{i=1}^{N_3} 1_{B_n(x_i , 2 \delta)} S_m ) \wedge \omega \geq   N_3 \times  \frac{ e^{- h_\nu(f) n  -2 \gamma n}}{2^{\underline{d}_\nu -\gamma } e^{(- \chi_1 n -12 \gamma n )(\underline{d}_\nu -\gamma ) }} \times  \\
\frac{C(\alpha_0,\gamma_0)}{2}e^{(-\chi_1n-15\gamma n)(\bar{d}_S+\gamma) +2\chi_2 n -18\gamma n} .
\end{align*}
So, as in the proof of Theorem \ref{Theorem1}: 
\begin{align*}
 C\lambda_1(f^n) \geq \lim_m\int (f^n)_*( S_m ) \wedge \omega \geq   \frac{e^{ h_\nu(f) n - 3\gamma n}}{8} \times \frac{ e^{- h_\nu(f) n  -2 \gamma n}}{2^{\underline{d}_\nu -\gamma }   e^{(- \chi_1 n -12 \gamma n )(\underline{d}_\nu -\gamma ) }}  \times  \\
\times  \frac{C(\alpha_0,\gamma_0)}{2}e^{(-\chi_1 n-15\gamma n)(\bar{d}_S+\gamma) +2\chi_2 n -18\gamma n}.
\end{align*}
Taking the logarithm and dividing by $n$ gives:
\begin{align*}
\frac{\log \lambda_1(f^n)}{n} \geq   \frac{1}{n} \log\left( \frac{C(\alpha_0, \gamma_0)}{C 16. 2^{\underline{d}_\nu -\gamma} }   \right)-23 \gamma + (\chi_1 + 12 \gamma)(\overline{d}_\nu- \gamma) -(\chi_1 + 15 \gamma)(\bar{d}_S+ \gamma) +2 \chi_2.
\end{align*}
We let $n \to \infty$ then $\gamma \to 0 $ then:
\begin{align*}
\log d_1 \geq   \chi_1\underline{d}_\nu- \chi_1\bar{d}_S +2 \chi_2
\end{align*}
which is what we wanted.

\noindent Henry De Th\'elin, 	Université Paris 13, \\
Sorbonne Paris Cité, LAGA \\
CNRS (UMR 7539), F-93430 \\ 
Villetaneuse, France \\
\noindent Email: dethelin@math.univ-paris13.fr  \\
  
\noindent Gabriel Vigny, LAMFA - UMR 7352, \\ 
U. P. J. V. 33, rue Saint-Leu, 80039 Amiens, France. \\
\noindent Email: gabriel.vigny@u-picardie.fr

\end{document}